\documentclass[11pt]{article}

\usepackage{graphicx,enumitem}
\usepackage{graphicx,psfrag,epsf}
{%
\centering\addtocounter{figure}{1}
\begin{enumerate}[%
itemsep=2pt,parsep=0em,
label={(\alph*)},
ref={\thefigure.(\alph*)}
]}%
{\end{enumerate}\addtocounter{figure}{-1}}

\usepackage[colorlinks=true,linkcolor=blue]{hyperref}%
\usepackage{graphicx}
\hypersetup{urlcolor=blue, citecolor=blue, colorlinks=true, linkcolor=blue} 
 
\usepackage[margin=1in]{geometry}

\usepackage[utf8]{inputenc}
\usepackage{amsmath}
\usepackage[english]{babel}
\usepackage{amssymb}
\usepackage{amsfonts}
\usepackage{bm}
\usepackage{graphicx}
\usepackage{amsmath}
\usepackage{amsthm}
\usepackage{xcolor}
\usepackage{parskip}
\usepackage{multicol}
\usepackage{float}

\usepackage{verbatim}

\floatstyle{plaintop}
\restylefloat{table}
\usepackage[tableposition=top]{caption}

\usepackage{rotating}
\usepackage{subcaption}
\usepackage{multirow}

\usepackage{chngcntr}
\usepackage{apptools}
\AtAppendix{\counterwithin{lem}{subsection}}
\AtAppendix{\counterwithin{proposition}{subsection}}
\AtAppendix{\counterwithin{equation}{subsection}}

\usepackage{authblk}

 \theoremstyle{definition}

\newtheorem{lem}{Lemma}[section]
\newtheorem{proposition}{Proposition}[section]

\newtheorem{prueba}{Proof}[section]
\newtheorem{ejemplo}{Example}[section]
\newtheorem{remark}{Remark}[section]

\usepackage{booktabs}
\usepackage{animate}

\usepackage{natbib}
\hypersetup{urlcolor=blue, citecolor=black, colorlinks=true, linkcolor=blue} 

\numberwithin{equation}{section}
\numberwithin{figure}{section}
\numberwithin{table}{section}

\usepackage{blindtext}
\usepackage{lineno}

\begin{document}

\def\spacingset#1{\renewcommand{\baselinestretch}%
{#1}\small\normalsize} \spacingset{1}

\title{{\bf Hybrid Parametric Classes of Isotropic Covariance Functions for Spatial Random Fields}}

\author[1]{Alfredo Alegr\'ia\footnote{Corresponding author. Email: alfredo.alegria@usm.cl}}

\author[1]{Fabián Ramírez}

\author[2]{Emilio Porcu}

\affil[1]{Departamento de Matem\'atica, Universidad T{\'e}cnica Federico Santa Mar{\'i}a, Chile}

\affil[2]{Department of Mathematics, 
Khalifa University, 
 The Arab Emirates}

\maketitle

\begin{abstract}
\noindent Covariance functions are the core of spatial statistics, stochastic processes, machine learning as well as many other theoretical and applied disciplines. The properties of the covariance function at small and large distances determine the geometric attributes of the associated Gaussian random field. Having covariance functions that allow to specify both local and global properties is certainly on demand. This paper provides a method to find new classes of covariance functions having such properties.  We term these models \emph{hybrid} as they are obtained as scale mixtures of piecewise covariance kernels against measures that are also defined as piecewise linear combination of parametric families of measures. In order to illustrate our methodology, we provide new families of covariance functions that are proved to be richer with respect to other well known families that have been proposed by earlier literature. More precisely, we derive a hybrid Cauchy-Mat{\'e}rn model, which allows us to index both long memory and mean square differentiability of the random field, and a hybrid Hole-Effect-Matérn model, which is capable of attaining negative values (hole effect), while preserving the local attributes of the traditional Matérn model. Our findings are illustrated through numerical studies with both simulated and real data.  \\


\end{abstract}

\section{Introduction}

Covariance functions are central to many disciplines such as spatial statistics \citep{Cressie:1993,Chiles2012,hristopulos2020random}, 
 stochastic processes \citep{porcu2018modeling, porcu2018shkarofsky}, 
 machine learning \citep{schaback2006kernel,james2013introduction, barp2022riemann},
   numerical analysis \citep{pazouki2011bases,cockayne2019bayesian} and stochastic mechanics \citep[][with the references therein]{ostoja2006material}. 
 Recent applications in climatology \citep{guinness2018compression,sacowea1}, oceanography \citep{furrer2007multivariate,horma2}, environmental sciences \citep{cressie2003spatial,stein2007spatial} and natural resources engineering \citep{chen2018current,emery2020geostatistics}  witness on the importance of covariance functions.

 It is very customary to assume the covariance function to depend on the distance between any pair of random variables located at two different points at the input space. Such an assumption is termed isotropy in spatial statistics and machine learning, and it is termed radial symmetry in other areas of applied mathematics. The behaviour of the covariance function at short or long distances (we call this local and global properties, respectively) is crucial to understand the properties of random processes with a given covariance function. Specifically, the local properties are related to the fractal dimension as well as the geometric properties (e.g., mean square differentiability) of the associated random process, as well as to its sample paths. On the other hand, the global behaviour of the covariance function allows to characterize persistency or antipersistency (i.e., the long term behaviour) of the associated process. Another global behaviour of great interest is the so-called hole effect, which means that the covariance function could take negative values in a certain interval.

 Finding parametric families of isotropic covariance functions that allow to index both local and global behaviour is a major challenge that has been tackled to a very limited extent. The Mat{\'e}rn family has been the cornerstone in spatial statistics for over half a century now \citep{stein-book}. Its popularity is due to a parameter that controls the degree of mean square differentiability and fractal dimension  of the corresponding random field \citep{stein-book}. Recently, \cite{bevilacqua2022unifying} have shown that the Mat{\'e}rn class is a special case of a richer class of models that, additionally to indexing local properties, allow to switch between compact or global supports. In turn, compactly supported models lead to sparse covariance matrices \citep{furrer2006covariance,kaufman2008covariance} and this implies considerable computational gains in both estimation and prediction. Unfortunately, the Mat{\'e}rn class does not allow to index global behaviour of the associated random process. The Generalized Cauchy family \citep{gneiting2004stochastic} allows to index the fractal dimension  and the long memory behaviour. Notably, it does not allow to index mean square differentiability, as the model is either non differentiable or infinitely differentiable at the origin. The same properties are shared by the Dagum model \citep{berg2008dagum}, which does not allow to index mean square differentiability either. None of the aforementioned models allow to attain negative spatial dependencies.

Spectral approaches can be a promising avenue to find flexible families of covariance functions.  \cite{laga2017modified} proposed a modified version of the spectral density associated with the Mat{\'e}rn family. The new class has two additional parameters that can be loosely interpreted as a continuous version of a moving average process. More recently, \cite{ma2022beyond} have proved that a two-fold application of Gaussian scale mixtures can provide models with polynomial decays while preserving the local properties of the candidate covariance function. Other non conventional properties of covariance functions have been studied by \cite{alegria2020crossdimple} and  \cite{alegria2021bivariate}, who proposed some modified scale mixtures representations to obtain classes of cross-covariance functions with non-monotonic behaviours (the so-called cross-dimple effect) for vector-valued random fields. In \cite{schlather2017parametric}, models that allow for a smooth transition between stationary and intrinsically stationary Gaussian random fields are derived.

All the previously mentioned parametric classes of covariance functions admit a scale mixture representation of a Gaussian kernel against a continuous, positive and bounded measure. Our paper starts from the Schoenberg integral representation of isotropic covariance functions on $\mathbb{R}^d$ \citep{schoenberg1938}, for all natural numbers, $d$. We specifically assume the Schoenberg measures to be parametric families of measures that are defined piecewise. Such a strategy is then shown to provide hybrid classes that generalize classes proposed in earlier literature. We illustrate this methodology by constructing a model that combines the global attributes of the Cauchy class and the local properties of the Matérn class. We show that the proposed model admits a closed form expression and examine its theoretical properties. Additionally, we study a more flexible formulation, where the Gaussian kernel involved in the scale mixture is replaced with a covariance kernel that is also  defined piecewise. Following this approach, we derive a hybrid model with local behaviour of Matérn type, and global behaviour that allows for covariance functions with negative values. We conduct numerical experiments with both simulated and real data in order to assess the statistical performance of the proposed models.

The article is organized as follows. Section \ref{background} contains a concise review of random fields and covariance functions coming from scale mixtures.  
Section \ref{sec3} presents general methodologies to build hybrid covariance models. Then, we derive the hybrid Cauchy-Matérn and the hybrid Hole-Effect-Matérn classes.
Section \ref{sec4} guides the reader through some numerical studies.  We finally provide a critical discussion in Section \ref{conclusions}, including a description of technical extensions of the present work such as the multivariate case, where covariance functions are matrix-valued, and the case of spherically indexed fields, where isotropy is defined in terms of the geodesic distance.


\section{Background}
\label{background}

\def\R{\mathbb{R}}
\def\S{\mathbb{S}}

Let $\{Z({\bf s}): {\bf s}\in \mathbb{R}^d\}$ be a (centered) second-order {stationary} Gaussian random field on $\mathbb{R}^d$. Such a field is completely characterized by its covariance function (or kernel). The isotropy of the covariance function is defined through a mapping $\varphi:[0,\infty) \rightarrow \mathbb{R}$ such that 
 $\text{cov}[ Z({\bf s}) ,Z({\bf s}')] = \varphi(h),$
  for every ${\bf s}, {\bf s}' \in\mathbb{R}^d$, where $h = \|{\bf s} - {\bf s}'\|$. The covariance function must satisfy the positive (semi) definiteness condition: for any $k\in\mathbb{N}$, $\{a_1,\hdots,a_k\}\subset\mathbb{R}$ and $\{{\bf s}_1,\hdots,{\bf s}_k\}\subset\mathbb{R}^d$,
  $$  \sum_{i,j=1}^k a_ia_j\varphi(\|{\bf s}_i-{\bf s}_j\|) \geq 0.$$
We use the notation $\varphi(\cdot; \bm{\lambda})$ for a parametric family of continuous covariance functions, where $\bm{\lambda}\in\mathbb{R}^p$ is a  vector of parameters. Further, we make use of the celebrated Schoenberg' theorem \citep{schoenberg1938}:  the functions $\varphi$ that are valid in any dimension $d\in\mathbb{N}$ are uniquely written as Gaussian scale mixtures of positive and bounded measures, that is  
\begin{equation*}
\label{mixture0}
\varphi(h; \bm{\lambda}) = \int_0^\infty  \exp(-uh^2) G({\rm d} u;\bm{\lambda}),  \qquad  h \ge 0,
\end{equation*}
where $\{ G({\rm d}\cdot; \bm{\lambda}), \; \bm{\lambda} \in \mathbb{R}^p \}$ is a parametric family of measures, that are termed Schoenberg measures in \cite{daley2014dimension}. Most of the covariance classes listed in the introduction admit such a representation against a measure that is absolutely continuous with respect to the Lebesgue measure, that is
\begin{equation}
\label{mixture}
\varphi(h; \bm{\lambda}) = \int_0^\infty  \exp(-uh^2) g(u;\bm{\lambda}) \text{d}u,   \qquad h \ge 0,
\end{equation}
for $\{ g(\cdot; \bm{\lambda}), \; \bm{\lambda} \in \mathbb{R}^p \}$ a parametric family of nonnegative functions. Throughout, we call $g$ the mixing function.

We now describe examples of some parametric classes of functions $\varphi$ that are determined according to (\ref{mixture}). Special attention is devoted to the Mat{\'e}rn, Cauchy and Generalized Cauchy models.  Other examples, including the stable and generalized hyperbolic models, can be found in \cite{yaglom1987correlation}, \cite{barndorff1978hyperbolic}, \cite{schlather2010some} and  \cite{porcu2018shkarofsky}.

\begin{ejemplo}[Matérn]
This class of covariance functions is defined as \citep{Matern}
\begin{equation}
\label{matern}
\varphi_{\mathcal{M}}\left(h ; \bm{\lambda} \right) =  \frac{2^{1-\nu}}{\Gamma(\nu)} (h/\alpha)^\nu K_\nu(h /\alpha), \qquad h \ge 0,
\end{equation}
where $\Gamma$ is the gamma function and $K_\nu$ is the modified Bessel function of the second kind \citep{Abramowitz-Stegun:1965}. Here, $\bm{\lambda}=[\alpha,\nu]^{\top}$, with
$\alpha$ and $\nu$ being positive parameters that control the scale (the rate of decay of the covariance in terms of $h$) and shape of (\ref{matern}), respectively. More precisely, $\nu$ regulates the degree of mean square differentiability of the random field (large values of $\nu$ are associated with smoother sample paths) \citep{stein-book}. When $\bm{\lambda}=[\alpha,1/2]^{\top}$, (\ref{matern}) simplifies into the exponential model, $\exp(- h / \alpha)$. On the other hand, as $\nu\rightarrow \infty$, a reparameterization of (\ref{matern}) tends to the Gaussian covariance function, defined as $\exp( - h^2/\alpha  )$.   
\end{ejemplo}

\begin{ejemplo}[Cauchy] This class of covariance functions is given by \citep{Chiles2012}
\begin{equation}
\label{cauchy}
\varphi_{\mathcal{C}}\left(h ; \bm{\lambda}\right) =    \left( 1 + h^2/\alpha  \right)^{-\nu/2}, \qquad h \ge 0,
\end{equation}
with $\bm{\lambda}=[\alpha,\nu]^{\top}$. As in the Matérn model, $\alpha>0$ is a scale parameter.  However, unlike the Matérn model which decays exponentially with distance \citep{stein-book}, (\ref{cauchy}) has a polynomial decay regulated by $\nu>0$. When $\nu\in(0,2)$, such a polynomial decay is connected with the Hurst parameter, a measure of long term memory, given by $H = 1-\nu/2$.  
\end{ejemplo}

 \begin{ejemplo}[Generalized Cauchy]
   This class of covariance functions is defined as (\citealp{gneiting2004stochastic} and the references therein)
\begin{equation}
\label{gcauchy}
\varphi_{\mathcal{GC}}\left(h ; \bm{\lambda}\right) =
\left( 1 +  h^\delta/\alpha  \right)^{-\nu/\delta}, \qquad h \ge 0,
\end{equation}
with $\bm{\lambda}=[\alpha,\nu,\delta]^{\top}$, where $\delta\in(0,2]$, $\alpha>0$ and  $\nu>0$. This generalized class preserves the polynomial decay of (\ref{cauchy}), but it is more flexible in the sense that the fractal dimension can be arbitrarily regulated through $\delta$ (see \citealp{gneiting2004stochastic} for details). Maybe surprisingly, this model does not allow to control the mean square differentiability of the respective random field, as the model is either non differentiable or infinitely differentiable at the origin.
\end{ejemplo}

Additional classes of covariance functions can be obtained from the more general mixture 
\begin{equation}
\label{mixture_gen}
\varphi(h; \bm{\lambda},\bm{\vartheta})= \int_0^\infty  \phi(h; u,\bm{\vartheta})  g(u;\bm{\lambda}) \text{d}u,   \qquad h \ge 0,
\end{equation}
where $\phi(\cdot; u,\bm{\vartheta})$ is an arbitrary covariance kernel, for every $u>0$, and $\bm{\vartheta}$ is a vector of parameters. Since the class of positive definite functions is a convex cone that is closed under the topology of pointwise convergence, if $\phi$ is valid (positive definite) in $\mathbb{R}^d$ for $d \leq d'$, for some $d'\in\mathbb{N}$, then $\varphi$ is valid in $\mathbb{R}^d$ for $d \leq d'$ as well. We refer the reader to  \cite{emery2006tbsim} for several explicit examples.


\section{Hybrid Classes of Covariance Functions}
\label{sec3}

\subsection{General Construction}
In this study, we propose new parametric classes of isotropic covariance functions,  $\widetilde{\varphi}(\cdot; \bm{\lambda},\bm{\omega},\bm{\xi})$,  determined according to  
\begin{equation}
    \label{expresion10}
   \widetilde{\varphi}(h; \bm{\lambda},\bm{\omega},\bm{\xi})  
  =  \omega_1  \int_0^{\xi_1} \exp(-uh^2) {g}_1(u; \bm{\lambda}_1) \text{d}u  + \omega_2 \int_{\xi_2}^\infty     \exp(-uh^2) {g}_2(u;\bm{\lambda}_2)  \text{d}u,
\end{equation}
where $g_1$ and $g_2$ are nonnegative functions on $[0,\xi_1)$ and $[\xi_2,\infty)$, respectively, and $\bm{\lambda} = [\bm{\lambda}_1^\top,\bm{\lambda}_2^\top]^\top$,  $\bm{\omega} = [\omega_1,\omega_2]^\top$ and $\bm{\xi} = [\xi_1,\xi_{2}]^\top$ are vectors of parameters, with $\omega_i, \xi_i>0$, for $i=1,2$. In other words, we replace the mixing function, $g$, in Equation (\ref{mixture}) with a function $\widetilde{g}$ that is defined piecewise, i.e.,
\begin{equation}
    \label{tuma_00}
     \widetilde{g}(u;\bm{\lambda},\bm{\omega},\bm{\xi})  = \omega_1   \, g_1(u;\bm{\lambda}_1) 1_{[0,\xi_1)}(u) +  \omega_2   \, g_2(u;\bm{\lambda}_2) 1_{[\xi_2,\infty)}(u), \qquad u \geq  0,
\end{equation}
with $1_A(\cdot)$ standing for the indicator function of a set $A$. Note that $\widetilde{g}$ may have discontinuities as it is built by gluing two individual pieces.  If the functions $g_i$ are continuous and bounded on their domains, a direct application of the dominated convergence theorem implies that the proposed covariance function (\ref{expresion10}) is continuous on $[0,\infty)$.  Throughout this manuscript, each function $g_i$ is positively proportional to a continuous probability density function.  Hence, the parametric family proposed in Equation (\ref{expresion10}) belongs to the Schoenberg class as defined through Equation (\ref{mixture}).

A more general construction considers different kernels in each segment of the mixture, i.e.,
\begin{equation}
    \label{expresion10_gen}
   \widetilde{\varphi}(h; \bm{\lambda},\bm{\omega},\bm{\xi},\bm{\vartheta})  
=  \omega_1  \int_0^{\xi_1} \phi_1(h; u, \bm{\vartheta}_1) {g}_1(u; \bm{\lambda}_1) \text{d}u  + \omega_2 \int_{\xi_2}^\infty    \phi_2(h; u, \bm{\vartheta}_2) {g}_2(u;\bm{\lambda}_2)  \text{d}u,    
\end{equation}
 where $\bm{\vartheta} =[\bm{\vartheta}^\top_1,\bm{\vartheta}^\top_2]^\top$. If $\phi_i$ is a valid  covariance function in $\mathbb{R}^d$ for $d\leq d'_i$, for some $d_i'\in\mathbb{N}$, $i=1,2$, then (\ref{expresion10_gen}) is a valid model in $\mathbb{R}^d$ if and only if  $d\leq \text{min}(d'_1,d'_2)$. The continuity of (\ref{expresion10_gen}) can be justified by following the same arguments used for the continuity of (\ref{expresion10}). 
 
\begin{remark}
Let us point out some additional remarks on this methodology.
\label{remark1}
\begin{enumerate}
\item When $\xi_1 = \xi_2 = \xi$, this parameter produces a continuous bridge  between two apparently disunited \emph{marginal} models. More precisely, as it increases from $0$ to $\infty$, we gradually go from $\omega_2 \int_{0}^\infty   \phi_2(h; u, \bm{\vartheta}_2) {g}_2(u;\bm{\lambda}_2)  \text{d}u$ to $\omega_1 \int_{0}^\infty   \phi_1(h; u, \bm{\vartheta}_1)  {g}_1(u;\bm{\lambda}_1)  \text{d}u$.
 \item 	When $\xi_1>\xi_2$, instead, there is a superposition of the {marginal} structures in the interval $[\xi_2,\xi_1)$. As $\xi_2\rightarrow 0$ and $\xi_1\rightarrow \infty$, we obtain the greatest possible superposition, which corresponds to a linear combination of the marginal models, $\omega_1 \int_{0}^\infty   \phi_1(h; u, \bm{\vartheta}_1) {g}_1(u;\bm{\lambda}_1)  \text{d}u + \omega_2 \int_{0}^\infty   \phi_2(h; u, \bm{\vartheta}_2) {g}_2(u;\bm{\lambda}_2)  \text{d}u.$
 \end{enumerate}
\end{remark}

The apparent flexibility of the proposed mixtures is justified by classical theory on local and global behaviour of covariance functions. In particular, a direct application of Tauberian theorems \citep{stein-book} proves that mean square differentiability of $\widetilde{\varphi}$ will be determined by $g_2$. On the other hand, direct inspection in concert with Equation (4) in \cite{gneiting2004stochastic}  shows that the long term behaviour of $\widetilde{\varphi}$ is decided by $g_1$.  The next sections show that it is possible to provide examples in algebraically closed form that allow to attain the desired flexibility.

\subsection{A Hybrid Cauchy-Matérn Class}
\label{closed-models}

We present a hybrid Cauchy-Mat{\'e}rn model, for which the acronym ${\cal CM}$ is used. This model is a special case of (\ref{expresion10}).  Let us first introduce the generalized incomplete gamma function  \citep{chaudhry1994generalized}, 
$$\Gamma(a;b;c)=\int_{b}^{\infty} t^{a-1} \exp(-t-c t^{-1}) \, \text{d}t,$$
and the lower incomplete gamma function, $\gamma(a,b) = \Gamma(a;0;0)-\Gamma(a;b;0)$. 
\begin{proposition} 
\label{botellazo}
Let  $\bm{\lambda} = [\bm{\lambda}_1^\top,\bm{\lambda}_2^\top]^\top$, with $\bm{\lambda}_i=[\alpha_i,\nu_i]^\top$,  $\bm{\omega} = [\omega_1,\omega_2]^\top$ and $\bm{\xi} = [{\xi}_1,\xi_2]^\top$ be vectors having positive elements. Let 
\begin{equation}
\label{closed_cm}
\widetilde{\varphi}_{\mathcal{CM}}(h; \bm{\lambda},\bm{\omega},\bm{\xi}) =  \omega_1  \, \widetilde{\varphi}^{\, (1)}_{\mathcal{C}}(h; \bm{\lambda}_1,\xi_1)    +  \omega_2  \,  \widetilde{\varphi}^{\, (2)}_{\mathcal{M}}(h; \bm{\lambda}_2,\xi_2),   \qquad h \ge 0,
\end{equation}
where 
\begin{equation}
\label{cauchy_p1}
\widetilde{\varphi}^{\, (1)}_{\mathcal{C}}(h; \bm{\lambda}_1,\xi_1)  =    \frac{\gamma(\nu_1/2,(h^2+\alpha_1)\xi_1)}{\Gamma(\nu_1/2)}   \varphi_{\mathcal{C}}(h; \bm{\lambda}_1) 
\end{equation} and 
\begin{equation}
\label{matern_p2}
\widetilde{\varphi}^{\, (2)}_{\mathcal{M}}(h; \bm{\lambda}_2,{\xi}_2)  =  \varphi_{\mathcal{M}}(h; \bm{\lambda}_2)  -   \frac{1}{\Gamma(\nu_2)} \Gamma\left( \nu_2; \frac{1}{4\xi_2 \alpha_2^2} ; \frac{h^2}{4 \alpha_2^2}\right),
\end{equation}
with $\varphi_{{\cal M}}$ and $\varphi_{{\cal C}}$ being, respectively, the Mat{\'e}rn and the Cauchy models defined at (\ref{matern}) and (\ref{cauchy}). Then, $\widetilde{\varphi}_{\mathcal{CM}}$ is positive definite in $\mathbb{R}^d$ for all $d\in\mathbb{N}$.
\end{proposition}

\begin{prueba}
We provide a proof of the constructive type, by showing that $\widetilde{\varphi}_{{\cal CM}}$ admits the representation (\ref{expresion10}), with $g_1(u;\bm{\lambda}_1) = g_{{\cal C}}(u;\bm{\lambda}_1)$ and $g_2(u;\bm{\lambda}_2) = g_{{\cal M}}(u;\bm{\lambda}_2)$, with $g_{{\cal C}}$ and $g_{{\cal M}}$ that are respectively defined as 
\begin{equation}
     \label{gamma}
g_{\mathcal{C}}(u;\bm{\lambda}_1) =  \frac{\alpha_1^{\nu_1/2}}{\Gamma(\nu_1/2)}  u^{\nu_1/2-1}\exp\left(-\alpha_1 u\right), 
\end{equation} and
\begin{equation}
     \label{inverse_gamma}
g_{\mathcal{M}}(u;\bm{\lambda}_2) = \frac{1}{\Gamma(\nu_2)} \left( \frac{1}{2 \alpha_2} \right)^{2\nu_2} u^{-\nu_2-1} \exp\left(-\frac{1}{4u \alpha_2^2}\right),
\end{equation} 
 where for both cases all the parameters are positive.  To attain the analytical expression of $\widetilde{\varphi}_{{\cal C}}^{(1)}$, we notice that 
\begin{eqnarray*}
\label{identity}
\int_{0}^{\xi_1}   \exp(- u h^2) {g}_{\mathcal{C}}(u;\bm{\lambda}_1) \text{d}u &=&  \varphi_{\mathcal{C}}(h; \bm{\lambda}_1)  \int_{0}^{\xi_1}  \frac{(h^2+\alpha_1)^{\nu_1/2}}{\Gamma(\nu_1/2)}  u^{\nu_1/2-1}\exp(-(h^2+\alpha_1) u) \text{d}u  \nonumber \\
&=& \varphi_{\mathcal{C}}(h; \bm{\lambda}_1)  \frac{\gamma(\nu_1/2,(h^2+\alpha_1)\xi_1)}{\Gamma(\nu_1/2)},
\end{eqnarray*}
where the second equality is due to the fact that the integral on the right hand side of the first line amounts to the cumulative distribution function of a gamma random variable with parameters $h^2+\alpha_1$ and $\nu_1/2$. 

To attain the expression of $\widetilde{\varphi}_{{\cal M}}^{(2)}$, we invoke Equation (10) in \cite{alegria2021bivariate}, so that
\begin{equation}
\label{identity2}
\int_{0}^{\xi_2}   \exp(- u h^2) {g}_{\mathcal{M}}(u;\bm{\lambda}_2) \text{d}u = \frac{1}{\Gamma(\nu_2)} \Gamma\left(\nu_2; \frac{1}{4\xi_2 \alpha_2^2}; \frac{h^2}{4\alpha_2^2}\right).
\end{equation}
The function $\widetilde{\varphi}_{{\cal M}}^{(2)}$ is thus attained by invoking formula 3.471.9 in  \cite{Grad}, for which we have
$\int_{0}^{\infty}   \exp(- u h^2) {g}_{\mathcal{M}}(u;\bm{\lambda}_2) \text{d}u = \varphi_{\mathcal{M}}(h;\bm{\lambda}_2).$
  \hfill $\qed$
\end{prueba}

When $\nu_2 = n + 1/2$, for some $n\in\mathbb{N}$,  (\ref{matern_p2}) can be expressed in terms of complementary error functions and modified Bessel functions of first and second kinds. We refer the reader to \cite{alegria2021bivariate} for a more detailed study of these special cases.

 The flexibility of the proposed structure is now illustrated through the following result, where we use the notation $f_1(h) \sim f_2(h)$, $h\rightarrow \infty$, to represent that, for some positive constant $c_0$, the  asymptotic relationship $\lim_{h\rightarrow \infty} f_1(h) / f_2(h)  = c_0$ holds.  
 
\begin{proposition}
 \label{botellazo2}
Let $Z$ be a Gaussian random field with covariance function of the form (\ref{closed_cm}). Then, $Z$ is $\kappa$-times mean square differentiable if and only if $\nu_2 > \kappa \ge 0$. Moreover, it is true that
$\widetilde{\varphi}_{\mathcal{CM}}(h; \bm{\lambda},\bm{\omega},\bm{\xi}) \sim h^{-\nu_1}$,  $h\rightarrow \infty.$
Hence, the Hurst parameter associated with $Z$ is solely indexed by the parameter $\nu_1$. 
\end{proposition}

 \begin{prueba} Arguments in Chapter 2 of  \cite{stein-book} show that an isotropic  random field with covariance function $\varphi$ is $\kappa$-times mean square differentiable if and only if $\varphi^{(2\kappa)}(0;\bm{\lambda})$ exists and is finite.
Direct inspection in concert with dominated convergence on the invoked Schoenberg representation (\ref{mixture}) show that this happens if and only if the mixing function $g$ satisfies
\begin{equation}
\label{msd_eq}
\int_0^\infty u^\kappa g(u;\bm{\lambda})\text{d}u < \infty.
 \end{equation} 
We use the latter argument for the special case of the function $\widetilde{\varphi}_{{\cal CM}}$, for which the tale of the resulting mixing function is uniquely determined by the mixing function associated with $\varphi_{{\cal M}}^{(2)}$ as in Proposition \ref{botellazo}. Direct inspection shows that (\ref{msd_eq}) is true  
 if and only if $\nu_2 > \kappa$. The first part of the proposition is established. 
 
For the second part, note that (\ref{cauchy_p1}) behaves as $h^{-\nu_1}$, as $h\rightarrow \infty$, because the lower incomplete gamma function involved in such an equation tends to $\Gamma(\nu_1/2)$, and the Cauchy class with parameter $\nu_1$ decays as $h^{-\nu_1}$. The result follows by noting that (\ref{matern_p2}) is dominated by the traditional Matérn model, which decays exponentially.
   \hfill $\qed$
\end{prueba}

To wrap up, the  hybrid Cauchy-Matérn model allows to index both mean square differentiability and long term behaviour of the associated Gaussian random field. We also note that these properties are independently addressed by the two parameters $\nu_1$ and $\nu_2$, and hence those parameters are statistically identifiable and allow to decouple local and global properties.

From a statistical viewpoint, a parsimonious choice may be considered  by setting  $\omega_1 = \omega_2=\omega$, $\alpha_1 = \alpha_2=\alpha$ and $\xi_1=\xi_2=\xi$. Thus, we obtain that Proposition \ref{botellazo} provides a five parameter family where $\omega$ indexes the variance, $\alpha$ the scale, $\nu_2$ the mean square differentiability, and $\nu_1$ the Hurst effect, whereas $\xi$ is a parameter that balances the shapes of the marginal structures involved in this model. Hence, (\ref{closed_cm}) generalizes the Mat{\'e}rn model in that it allows for polynomial decay while indexing continuously mean square differentiability.  

\begin{figure}
    \centering
    \includegraphics[scale=0.4]{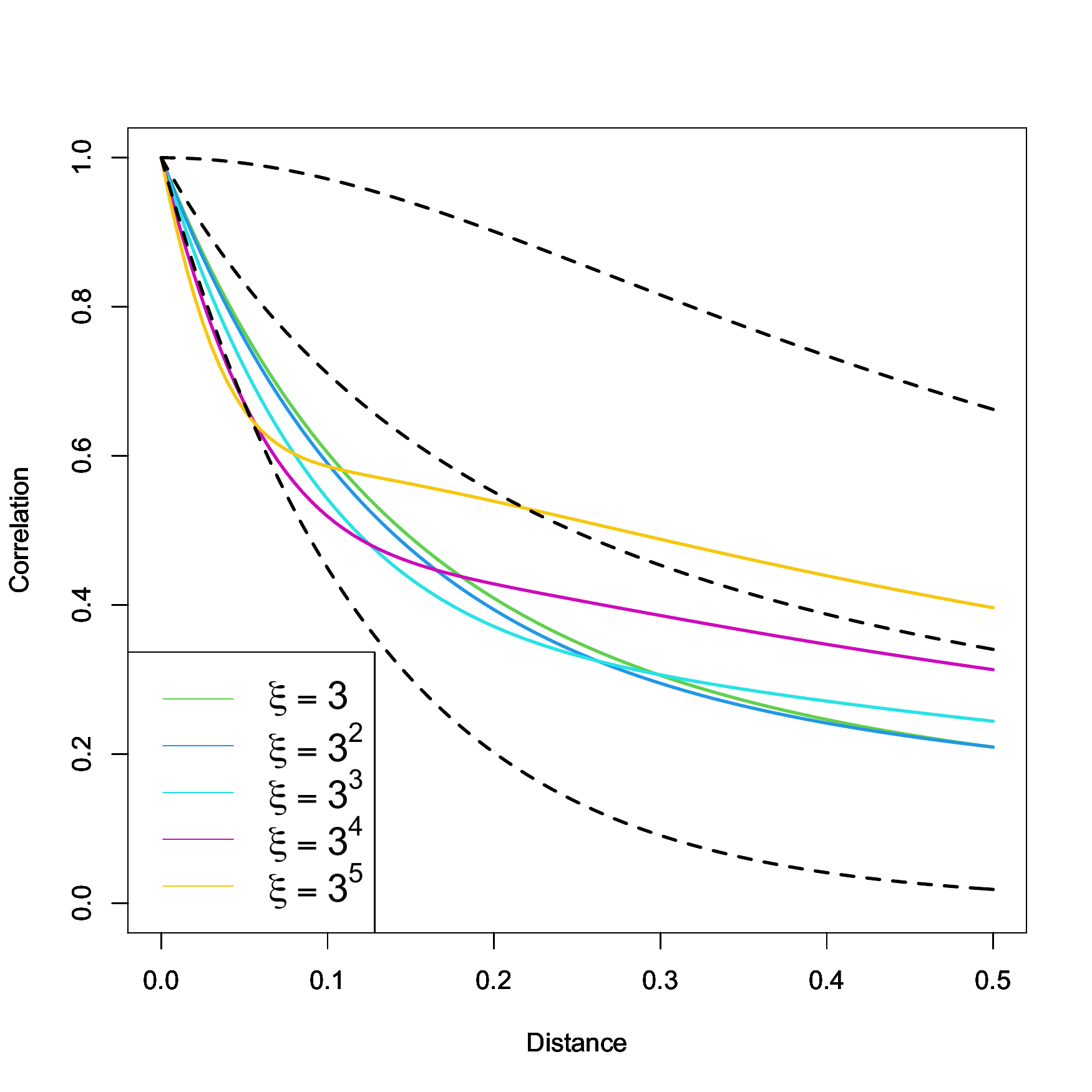}        \includegraphics[scale=0.4]{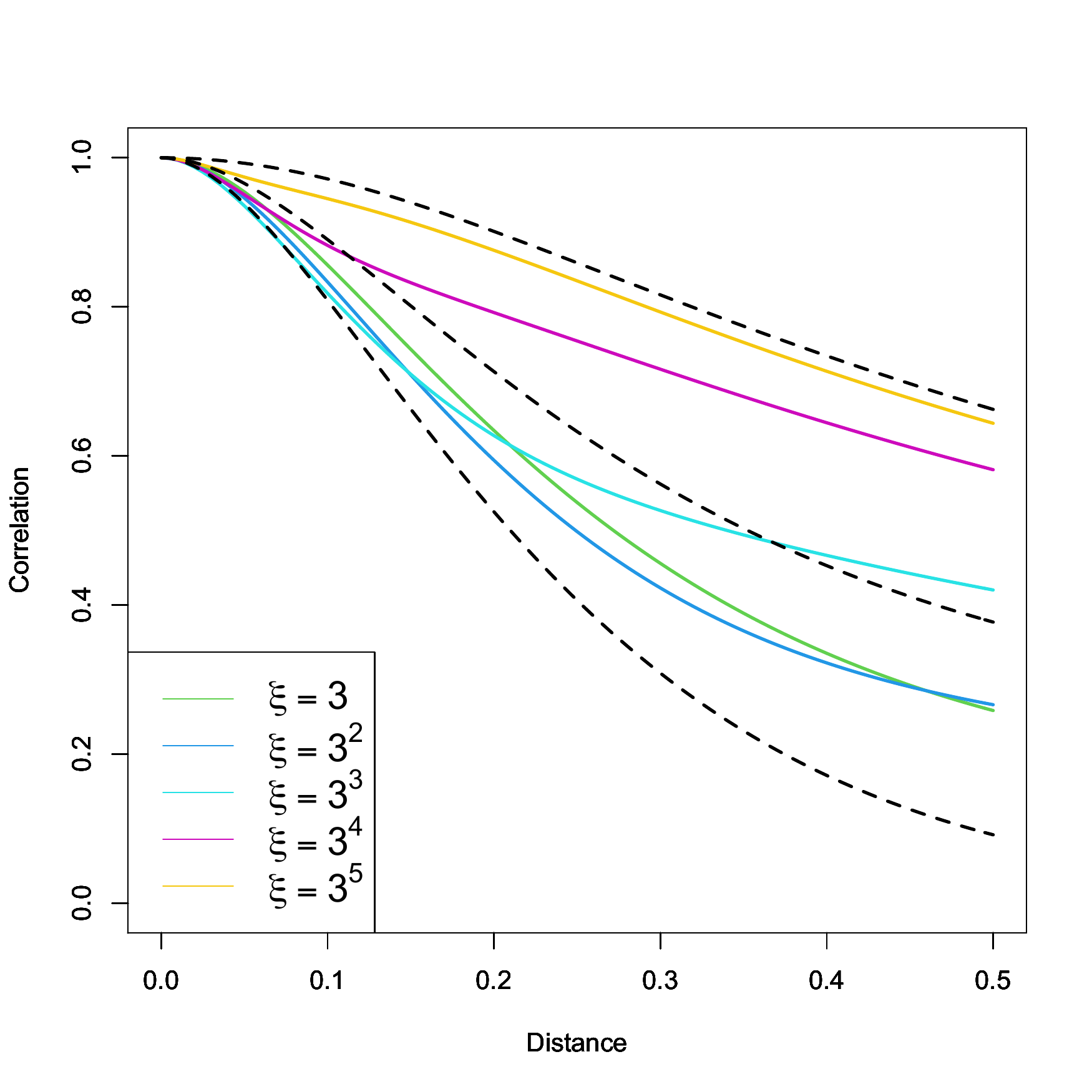}    
    \caption{Parsimonious hybrid Cauchy-Matérn model for $\omega=1/2$, $\alpha=1/8$,  $\nu_1=3/4$ and different values of $\xi$.  (Left)  $\nu_2=1/2$ and (Right)  $\nu_2=3/2$. The dashed lines represent the purely Cauchy, purely Matérn, and their average. All the models have been appropriately rescaled in order to obtain correlation functions.}
    \label{fig:curvas}
\end{figure}

Figure \ref{fig:curvas} shows the parsimonious hybrid Cauchy-Matérn model for different values of $\xi$. The traditional Matérn and traditional Cauchy, as well as their average, which are also special cases of the hybrid construction, are reported for comparison purposes. Note that the curves have a linear or parabolic decay near the origin according to $\nu_2=1/2$ or $\nu_2=3/2$, respectively, and then the decay is more gradual (polynomial rate) for large distances according to $\nu_1$, which is consistent with the local and global patterns that are coexisting.  We observe that $\xi$ has a manifest impact on the shape of the covariance function, as it produces some interesting forms (apparent changes of concavity) that could be useful in practice.

\subsection{A Hybrid Hole-Effect-Matérn Class}

We now present a hybrid class of covariance functions, with local attributes of Matérn type, attaining negative values at large distances. We use the acronym $\mathcal{HM}$ for this model, termed hybrid Hole-Effect-Matérn.  The proposed class comes from the mixture (\ref{expresion10_gen}), where $\phi_1$ is chosen in such a way that the resulting model can take negative values. 
\begin{proposition}
 \label{new_thm}
Let  $\bm{\lambda} = [\bm{\lambda}_1^\top,\bm{\lambda}_2^\top]^\top$, with $\bm{\lambda}_i=[\alpha_i,\nu_i]^\top$,  $\bm{\omega} = [\omega_1,\omega_2]^\top$ and $\bm{\xi} = [{\xi}_1,\xi_2]^\top$ be vectors having positive elements, and $\bm{\vartheta} = [\tau,\eta]^\top$ be a vector of additional parameters. Let 
\begin{equation}
\label{closed_hm}
\widetilde{\varphi}_{\mathcal{HM}}(h; \bm{\lambda},\bm{\omega},\bm{\xi},\bm{\vartheta}) =  \omega_1  \, \widetilde{\varphi}_{\mathcal{H}}^{\, (1)}(h; \bm{\lambda}_1,\xi_1,\bm{\vartheta})    +  \omega_2  \,  \widetilde{\varphi}_{\mathcal{M}}^{\, (2)}(h; \bm{\lambda}_2,\xi_2),   \qquad h \ge 0,
\end{equation}
where 
\begin{equation}
\label{h_p1}
\widetilde{\varphi}^{\, (1)}_{\mathcal{H}}(h; \bm{\lambda}_1,\xi_1,\bm{\vartheta})  =  \frac{\tau}{\Gamma(\nu_1)} \Gamma\left(\nu_1; \frac{1}{4\xi_1\alpha_1^2}; \frac{\eta   h^2}{4\alpha_1^2}\right)   - \frac{1}{\Gamma(\nu_1)} \Gamma\left(\nu_1; \frac{1}{4\xi_1 \alpha_1^2}; \frac{h^2}{4 \alpha_1^2}\right),
\end{equation}
and $\widetilde{\varphi}^{\, (2)}_{\mathcal{M}}$ as in (\ref{matern_p2}).  Then, $\widetilde{\varphi}_{\mathcal{HM}}$ is positive definite in $\mathbb{R}^d$ if and only if $1<\eta < \tau^{2/d}$.
\end{proposition}

\begin{prueba} We consider the construction (\ref{expresion10_gen}), with both $g_1$ and $g_2$ of the form (\ref{inverse_gamma}), and $\phi_2$ of Gaussian type. Thus, the derivation of $\widetilde{\varphi}^{\, (2)}_{\mathcal{M}}$ follows the same arguments employed in the proof of Proposition \ref{botellazo}. 

Before deriving (\ref{h_p1}), let us introduce the following lemma, which is a combination of Corollaries 4, 8 and 11 in \cite{posa}.  
\begin{lem}
\label{lema_posa}
The mapping $h\mapsto A \exp(-a h^2)  - B \exp(-b h^2)$ is positive definite in $\mathbb{R}^d$ if and only if 
\begin{equation}
\label{condition_posa}
1 < \frac{a}{b} <\left( \frac{A}{B} \right)^{2/d}.
\end{equation} 
\end{lem}

Although \cite{posa} focused on dimensions $d\leq 3$,  the same proof can be used in arbitrary dimensions.  To obtain the expression (\ref{h_p1}), we take the following covariance kernel in the first segment of the scale mixture
 \begin{equation}
 \label{neg_gaussian}
 \phi_1(h; u, \bm{\vartheta}) = \tau \exp(-u\eta h^2) - \exp(-uh^2), \qquad h\geq 0.
 \end{equation}
 Lemma \ref{lema_posa}  ensures that (\ref{neg_gaussian}) is positive definite in $\mathbb{R}^d$, provided that $u>0$ and $1 < \eta <\tau^{2/d}$. Thus,
 \begin{equation}
 \label{aa}
 \widetilde{\varphi}^{\, (1)}_{\mathcal{H}}(h; \bm{\lambda}_1,\xi_1,\bm{\vartheta})= \tau  \int_0^{\xi_1}  \exp(-u\eta h^2) g_{\mathcal{M}}(u;\bm{\lambda}_1) \text{d}u -  \int_0^{\xi_1}  \exp(-u h^2) g_{\mathcal{M}}(u;\bm{\lambda}_1) \text{d}u.
 \end{equation}
 Finally, we invoke the identity (\ref{identity2}), and we apply it to each integral involved in the right hand side of Equation  (\ref{aa}).
   \hfill $\qed$
\end{prueba}

The covariance function (\ref{neg_gaussian}) always takes negative values \citep{posa}, so it is a natural building block to achieve hybrid models with hole effect. The parameters in $\bm{\vartheta}$ are responsible for the sharpness of the hole effect. More precisely, as $\eta$ approaches $\tau^{2/d}$, the hole effect is more pronounced because the positive term in the right hand size of (\ref{neg_gaussian}) has less dominance. Moreover, when $d=1$, we have the least restrictive condition on $\eta$, and the resulting hole effect is more marked.  It is well known that the possibility of significant negative correlations vanishes as the dimension increases (see, e.g., page 45 in \citealp{stein-book}).

The next proposition characterizes the local attributes of (\ref{closed_hm}) and provides a lower bound for this model.

\begin{proposition}
\label{botellazo3}
Let $Z$ be a Gaussian random field with covariance function of the form (\ref{closed_hm}). Then, $Z$ is $\kappa$-times mean square differentiable if and only if $\nu_2 > \kappa \ge 0$. Moreover, we have the lower bound
\begin{equation}
\label{lower_bound}
\widetilde{\varphi}_{\mathcal{HM}}(h; \bm{\lambda},\bm{\omega},\bm{\xi},\bm{\vartheta}) \geq   \omega_1 (\tau \eta)^{-1/(\eta-1)} \left(\frac{1 - \eta}{\eta}\right) \left( 1 -  \frac{\gamma(\nu_1;\alpha_1/\xi_1)}{\Gamma(\nu_1)}   \right), \qquad h\geq 0.
\end{equation}
\end{proposition}

\begin{prueba}
The fact that $\nu_2$ controls the mean square differentiability is a direct consequence of the arguments used in the proof of Proposition \ref{botellazo2}. On the other hand, to find a lower bound, we note that 
\begin{eqnarray*}
\widetilde{\varphi}_{\mathcal{HM}}(h; \bm{\lambda},\bm{\omega},\bm{\xi},\bm{\vartheta})   &     \geq    &       \omega_1  \, \inf_{h\geq 0} \widetilde{\varphi}_{\mathcal{H}}^{\, (1)}(h; \bm{\lambda}_1,\xi_1,\bm{\vartheta})    +  \omega_2  \,  \inf_{h\geq 0}   \widetilde{\varphi}_{\mathcal{M}}^{\, (2)}(h; \bm{\lambda}_2,\xi_2)\\
& = & \omega_1 \int_0^{\xi_1}   \inf_{h\geq 0}  \phi_1(h; u, \bm{\vartheta})   g_1(u;\bm{\lambda}_1)\text{d}u.
\end{eqnarray*}
A straightforward calculation shows that $\phi_1$ attains its minimum value at $h^*= \sqrt{\frac{\log(\tau \eta)}{u(\eta-1)}}$. Thus, 
$$ \phi_1(h; u, \bm{\vartheta}) \geq   \phi_1(h^*; u, \bm{\vartheta}) =  \tau \exp\left(- \frac{\eta  \log(\tau \eta)}{\eta-1}\right) - \exp\left(- \frac{\log(\tau \eta)}{\eta-1} \right) =  (\tau \eta)^{-1/(\eta-1)} \left(\frac{1 - \eta}{\eta}\right).$$
Since $g_1$ is given by (\ref{inverse_gamma}), we invoke the formula of the cumulative distribution function of an inverse gamma random variable to establish that
$$ \int_{0}^{\xi_1} g_{1}(u;\bm{\lambda}_1) \text{d}u = 1 -  \frac{\gamma(\nu_1;\alpha_1/\xi_1)}{\Gamma(\nu_1)}. $$
The proof is completed.
  \hfill $\qed$
\end{prueba}

Note that as $\xi_1\rightarrow \infty$ (i.e., as the hole effect predominates), the lower bound in Equation (\ref{lower_bound}) decreases to $(\tau \eta)^{-1/(\eta-1)} (1-\eta)/\eta$. On the contrary, as $\xi_1\rightarrow 0$, such a bound increases to zero, i.e., the hole effect becomes negligible, which is not surprising, because in such a case the Matérn class is predominant. A similar conclusion can be obtained in the limit case $\eta\rightarrow 1$.

 A parsimonious variant of this model consists of taking  $\omega_1 = \omega_2=\omega$ (variance parameter), $\alpha_1=\alpha_2=\alpha$ (scale parameter) and $\nu_1=\nu_2=\nu$ (smoothness parameter), whereas $\bm{\vartheta}$ regulates the hole effect (as discussed above) and  $\xi_1=\xi_2=\xi$ has a similar interpretation as in the hybrid Cauchy-Matérn model.

\begin{figure}
    \centering
    \includegraphics[scale=0.4]{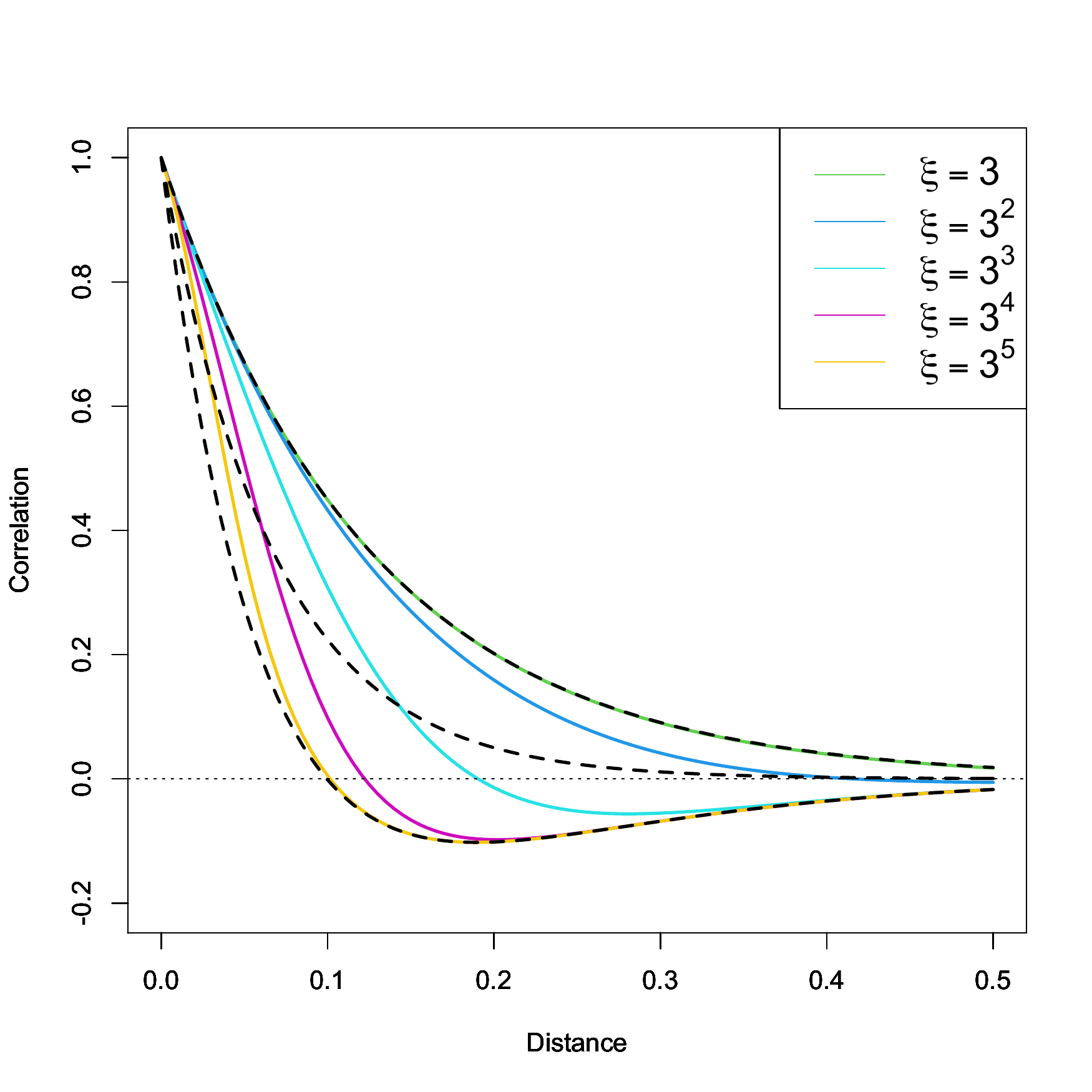}        \includegraphics[scale=0.4]{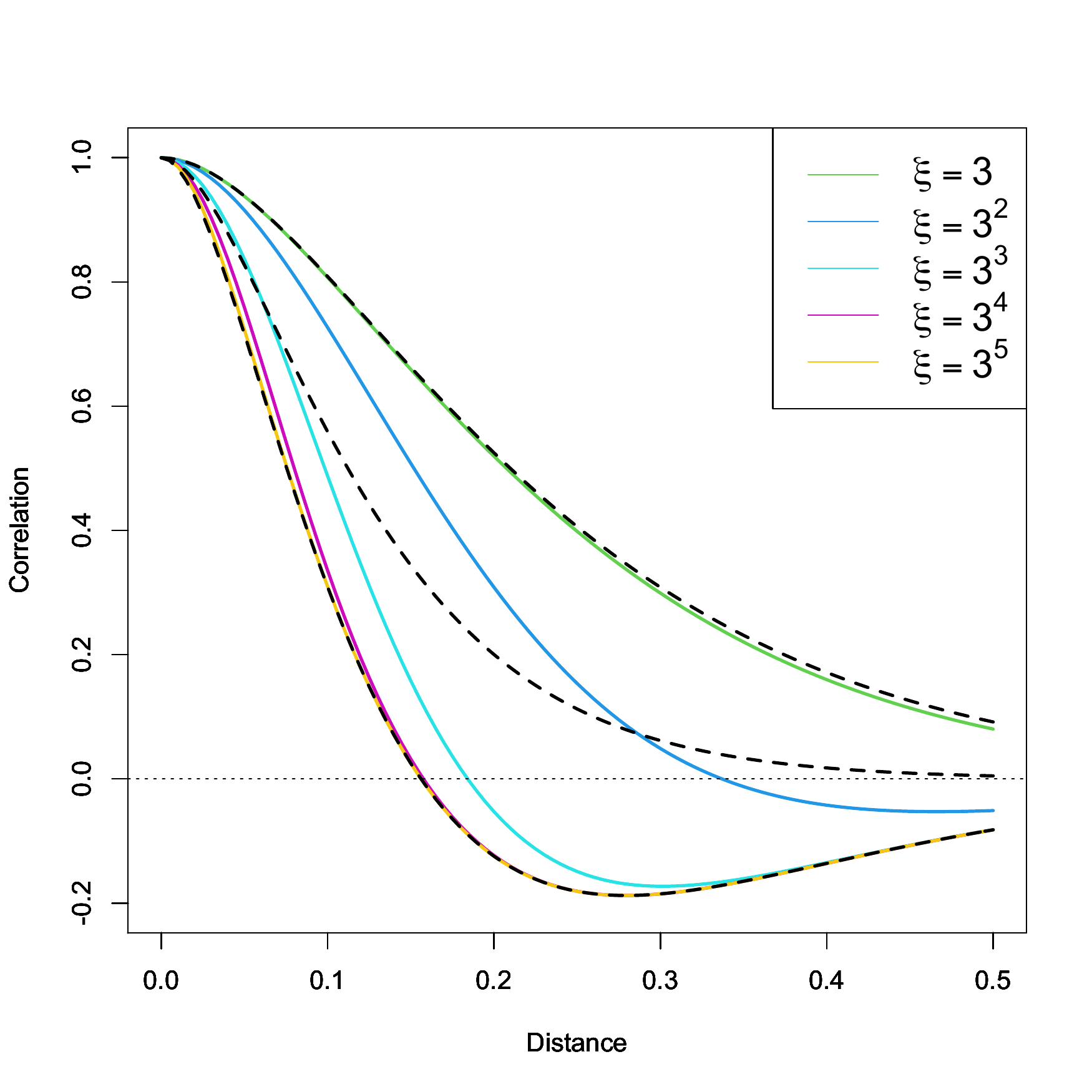}    
    \caption{Parsimonious hybrid Hole-Effect-Matérn  model in dimension one, for $\omega=1/2$, $\alpha=1/8$, $\tau = 2$, $\eta=7/2$ and different values of $\xi$.  (Left)  $\nu=1/2$ and (Right)  $\nu=3/2$. The dashed lines represent the limit cases reported in Remark \ref{remark1}. All the models have been appropriately rescaled in order to obtain correlation functions.}
    \label{fig:curvas2}
\end{figure}

Figure \ref{fig:curvas2} shows the parsimonious hybrid Hole-Effect-Matérn model for different values of $\xi$. The limit cases described in Remark \ref{remark1} are also reported, in a similar fashion to Figure \ref{fig:curvas}. It can be seen that negative values coexist with different levels of smoothness at the origin, as expected.

\section{Numerical Experiments}
\label{sec4}

\subsection{Simulated Data}
\label{simula}

We conduct simulation studies to assess the performance of maximum likelihood inference when a hybrid covariance structure is present. We focus on the parsimonious hybrid Cauchy-Matérn dependence structure, as it will be applied to real data in the next section. We consider $\omega = 1$, $\alpha=1/8$, $\nu_1=3/4$ and the following scenarios for $[\nu_2,\xi]$:  (a) $[1/2,40]$, (b) $[1/2,120]$, (c) $[3/2,40]$ and (d) $[3/2,120]$. For each scenario, we simulate 200 independent realizations of a Gaussian random field on 100 uniformly sampled points  in the square $[0,3]^2$ and estimate the parameters through maximum likelihood. We then repeat the experiment with 256 spatial locations. We only estimate $\omega$ , $\alpha$ and $\xi$, whereas $\nu_1$ and $\nu_2$ are fixed, which is a common practice in geostatistics. Instead of directly estimating $\xi$, we consider the following alternative parameterization: $\widetilde{\xi} = \sqrt{\xi}  \alpha$, which seems to be a  natural choice according to Equations (\ref{cauchy_p1}) and (\ref{matern_p2}). 

\begin{figure}
\centering
\includegraphics[scale=0.35]{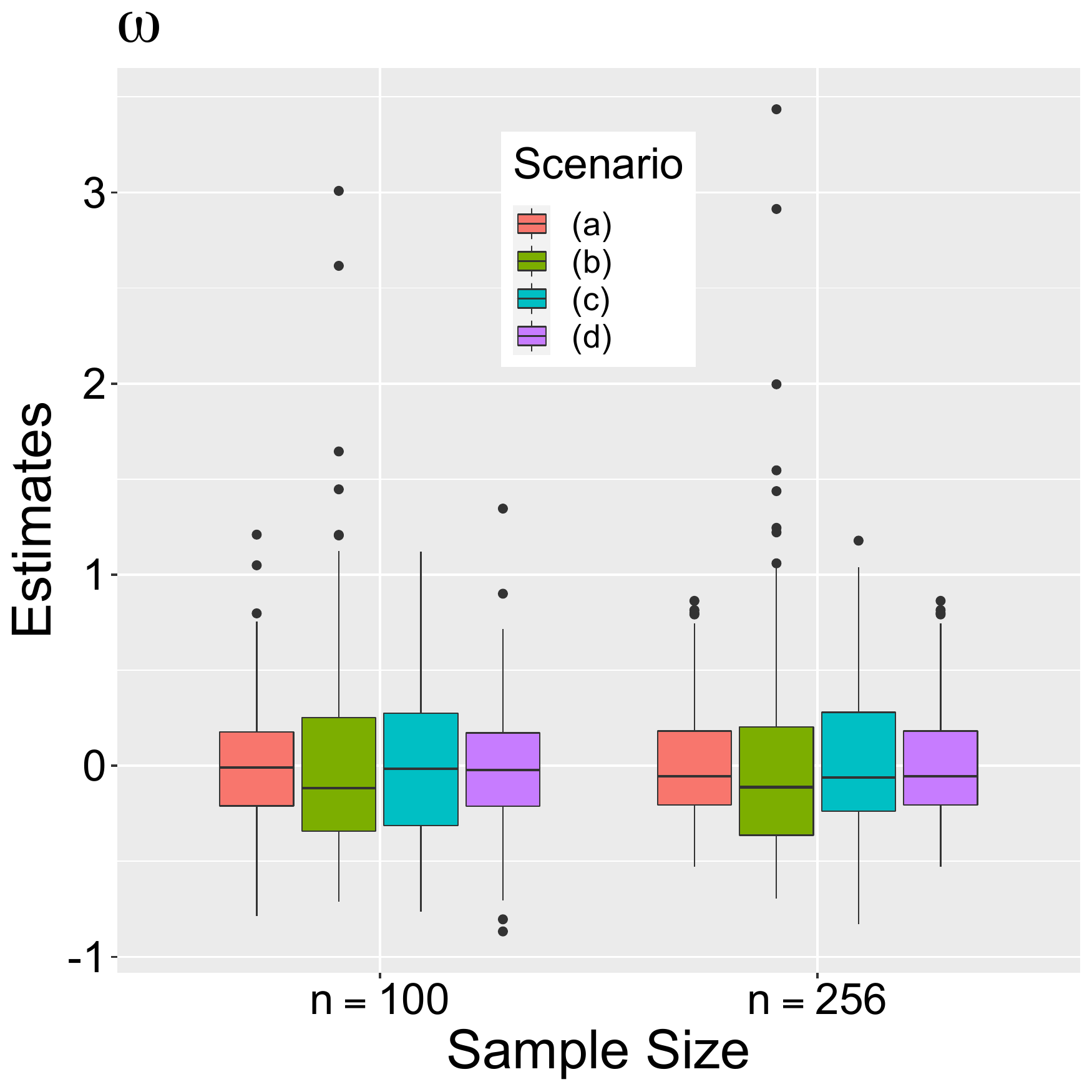}
\includegraphics[scale=0.35]{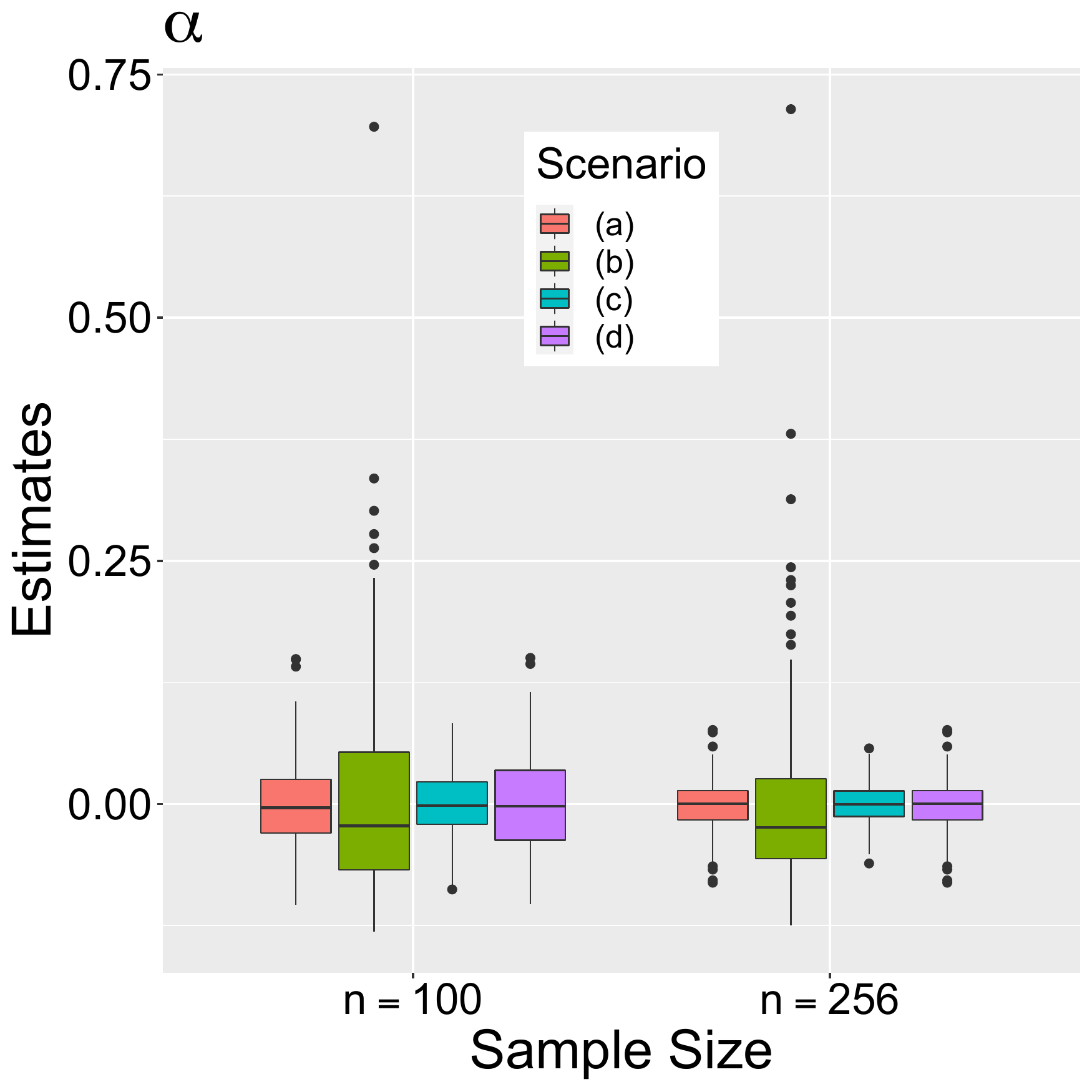}
\includegraphics[scale=0.36]{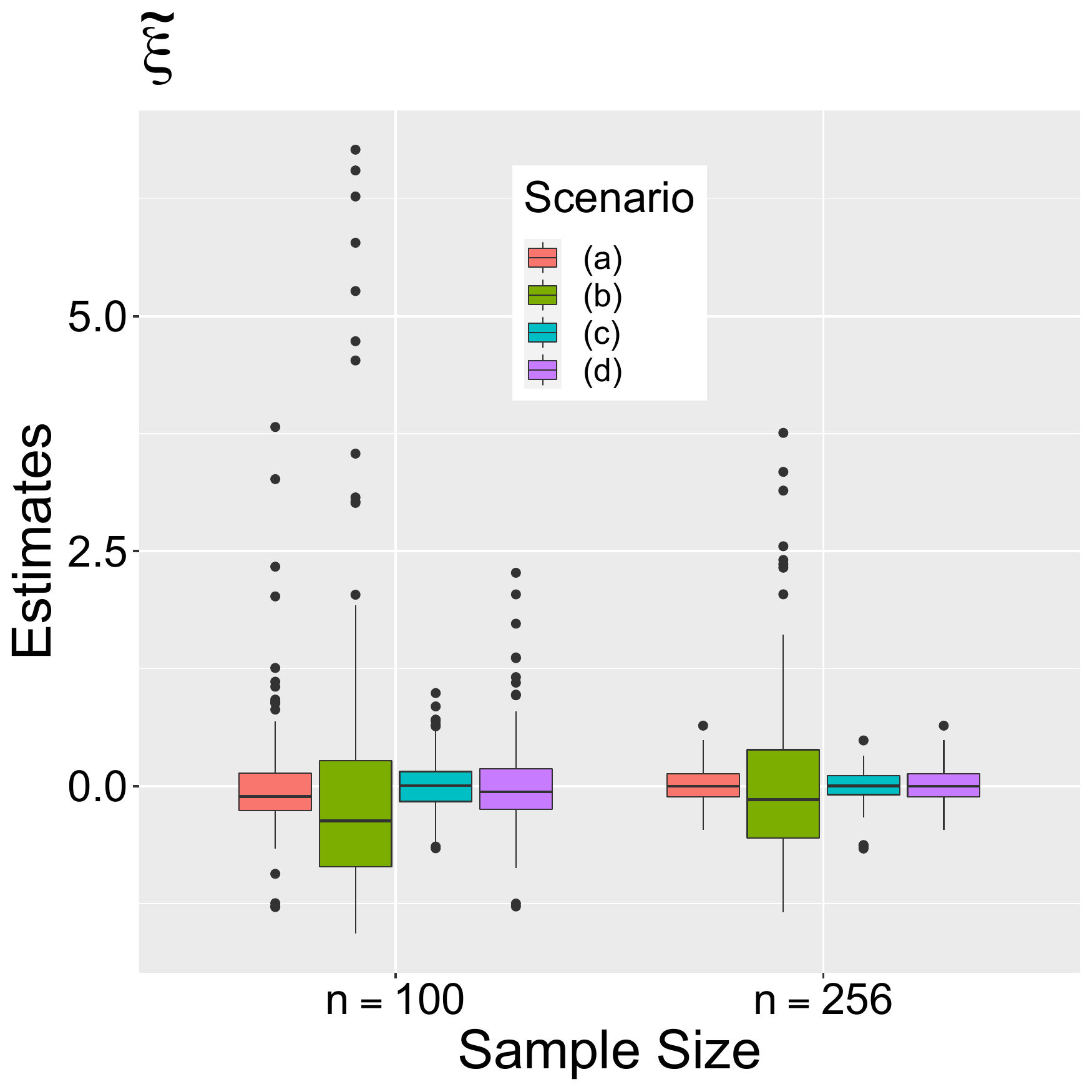}
\caption{Centered boxplots of the maximum likelihood estimates for  the parsimonious hybrid Cauchy-Matérn model in scenarios (a)-(d).}
\label{boxplots}
\end{figure}

\begin{figure}
\centering
\includegraphics[scale=0.4]{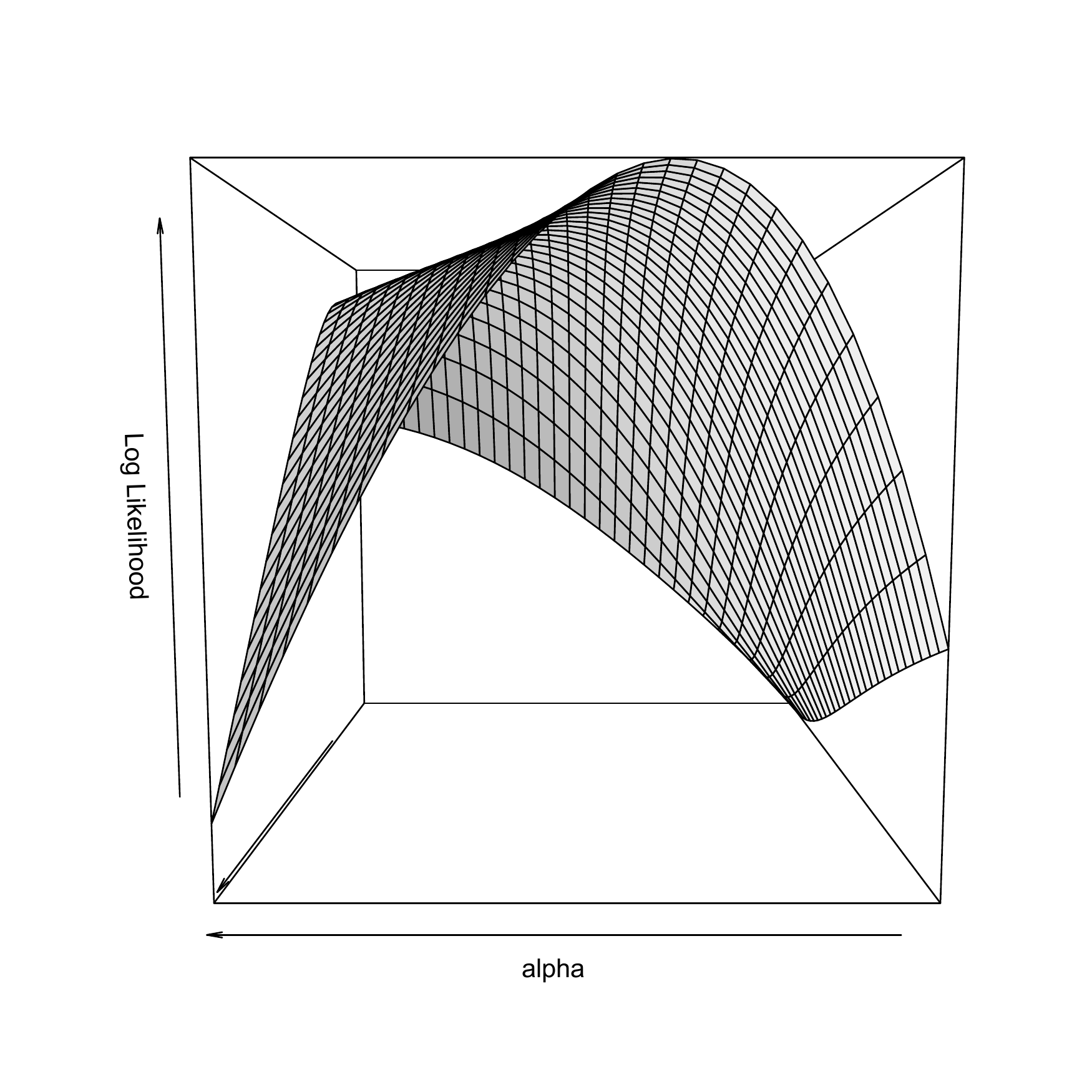}
\includegraphics[scale=0.4]{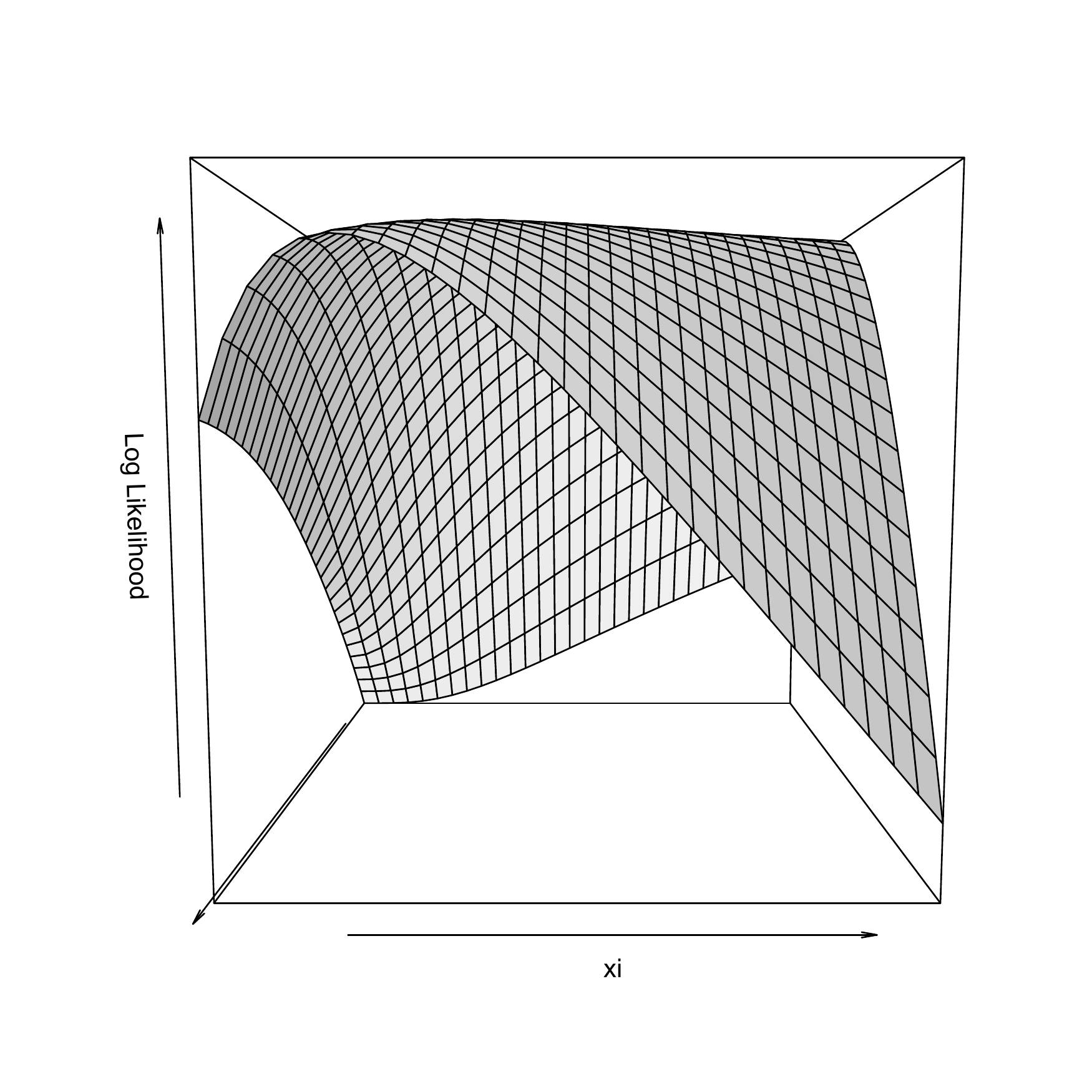}
\caption{Log-likelihood function, with respect to $\alpha$ and $\xi$, for scenario (b).  Left and right panels correspond to the same plot from different viewpoints. }
\label{loglik_profile}
\end{figure}

Figure \ref{boxplots} displays the results. The estimates are approximately unbiased and the variance decreases as the sample size increases from 100 to 256, which is an expected behaviour.  The variability of the estimates substantially decreases in scenarios (c) and (d), i.e., when the random field is smoother, which is a typical attribute of likelihood-based estimates in this context  \citep{bevilacqua2015comparing}.  On the contrary, such a  variability deteriorates as $\xi$ increases from $40$ to $120$. Figure \ref{loglik_profile} shows the log-likelihood in terms of $\xi$ and $\alpha$, with fixed $\omega$, for a single realization of the random field, under scenario (b). Although the surface has a clear maximum value, the objective function is apparently more flat in the direction of $\xi$. This could explain the increased variability  in scenarios (b) and (d), with respect to (a) and (c). Despite the previous remarks, in general, the estimates appear to be reasonable in each scenario and no identifiability issues are observed.

We now explore the predictive performance of the proposed class through a cross validation analysis. We simulate 200 independent realizations on 100 uniformly sampled locations in $[0,3]^2$ according to the scenarios (a)-(d) described above. We assess the accuracy through a leave-one-out prediction strategy in terms of the mean squared error (MSE), mean absolute error (MAE),  log-score (LSCORE) and continuous ranked probability score (CRPS) (see \citealp{zhang2010kriging}). Small values of these indicators suggest superior predictions. We evaluate the performance of the hybrid Cauchy-Matérn model, using the Generalized Cauchy class as benchmark. Thus, for each realization, we estimate the parameters with both models and proceed to make the predictions through a simple kriging approach.  The Generalized Cauchy model (\ref{gcauchy}) has been augmented with a multiplicative parameter $\omega$, namely $h \mapsto \omega (1+ h^\delta/\alpha)^{-\nu/\delta}$, so it is parameterized by $\omega$ and ${\alpha}$, and $\nu=3/4$ and $\delta=1,2$ are fixed. 

Table \ref{pred_tab} shows that, in each scenario, the proposed hybrid model outperforms its competitor.  All the cross-validation scores substantially decrease in scenarios (c) and (d). From this brief study, we observe  that when the true underlying covariance has a hybrid  structure, an incorrect specification of the spatial association has a negative impact on the posterior predictions. Since the behaviour of an isotropic covariance function near the origin has a strong impact on the quality of predictions \citep{stein-book}, our simulation experiment suggests that in some circumstances the local shape of the proposed model cannot be replicated by other appealing existing structures.

\begin{table}
\renewcommand{\arraystretch}{1.3}
\centering
\caption{Cross-validation scores for the parsimonious hybrid Cauchy-Matérn  and Generalized Cauchy (with $\delta=1,2$) models in scenarios (a)-(d).}
\begin{tabular}{cccccc}\hline \hline
    Scenario  &  Model                     &  MSE   &  MAE   & LSCORE & CRPS\\ \hline
               (a)  &     Hybrid Cauchy-Matérn &  0.706 &  0.668  & 1.231  &  1.773 \\               
                   &     Generalized Cauchy ($\delta=1$) &  0.714 & 0.672 &  1.238  &  1.787    \\
                   &     Generalized Cauchy ($\delta=2$)  &  0.718 & 0.674 & 1.241 & 1.798   \\
\hline
               (b)   &     Hybrid Cauchy-Matérn &  0.480 & 0.549 & 1.034 &  1.462   \\
               
                   &     Generalized Cauchy ($\delta=1$) &  0.489 & 0.555 & 1.046 & 1.480 \\
                   &     Generalized Cauchy ($\delta=2$)  &  0.497 & 0.559 &  1.055  & 1.506  \\ 
 \hline
               (c)   &     Hybrid Cauchy-Matérn &   0.172 &  0.316 & 0.398 & 0.840    \\
               
                   &     Generalized Cauchy ($\delta=1$) &    0.176  & 0.319 & 
 0.446 &  0.851   \\
                   &     Generalized Cauchy ($\delta=2$)  &  0.176 & 0.319 &
  0.414  &  0.862  \\ 
\hline
               (d)   &     Hybrid Cauchy-Matérn &    0.075 &  0.210 & 0.024 &  0.561   \\
               
                   &     Generalized Cauchy ($\delta=1$) &   0.077 &  0.213 & 
 0.052 &  0.575  \\
                   &     Generalized Cauchy ($\delta=2$)  &  0.082 & 0.219 &
 0.103  &   0.612   \\  \hline
\end{tabular}
\label{pred_tab}
\end{table}

\subsection{A Real Data Illustration}

  The estimation of recoverable resources is a task of fundamental importance in modern mining processes. A sound evaluation of such resources is crucial from an economic viewpoint and is critical for assessing the long-term availability of mineral resources and its impact on society. We consider a data set from a lateritic nickel deposit mined by open pit in Colombia, which contains measurements of the grades  of nickel, iron, chrome, alumina, magnesia and silica. 
  
\begin{figure}
\centering
\includegraphics[scale=0.6]{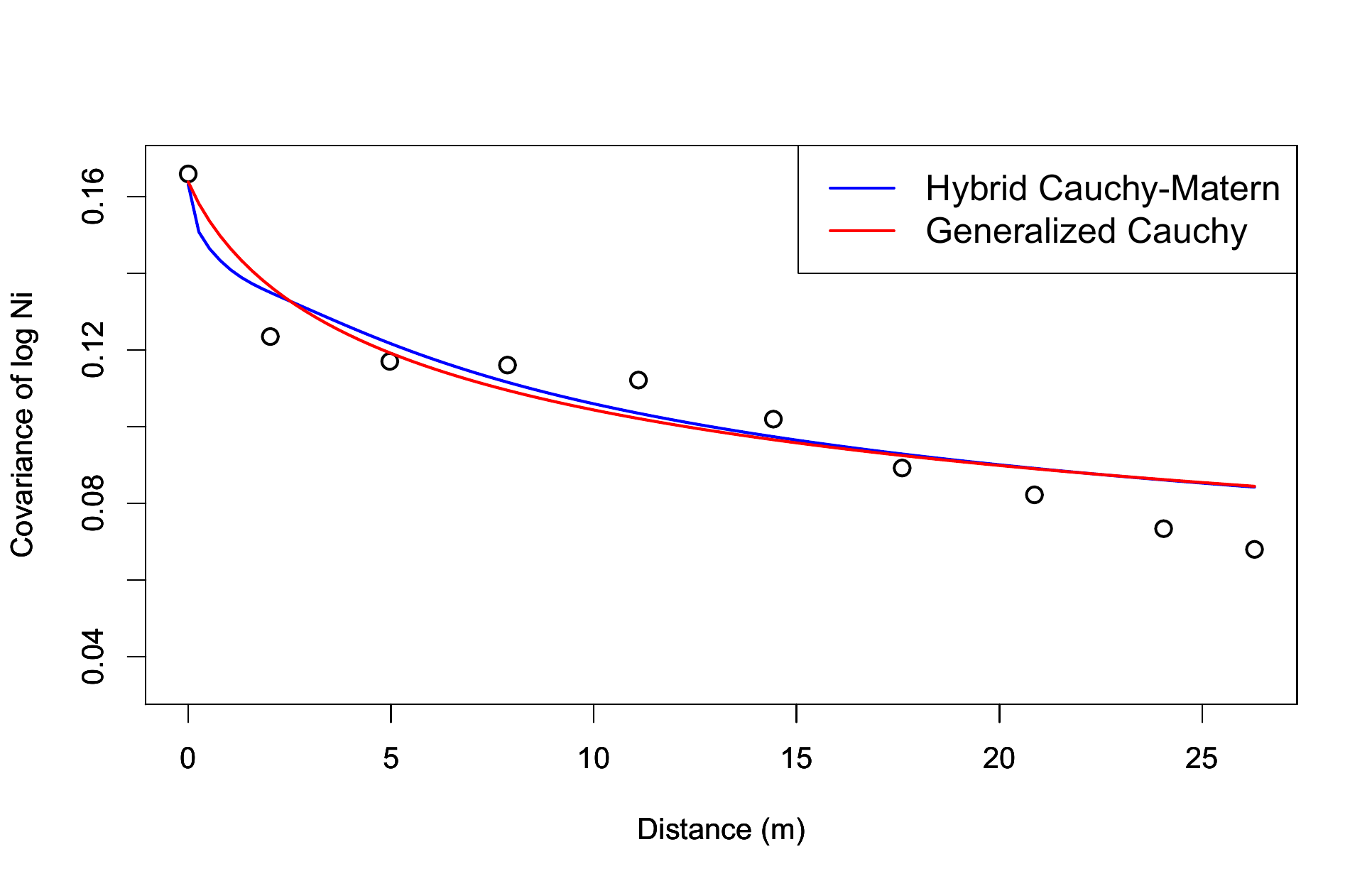}
\caption{Sample (circles) and modeled (solid lines) covariances of log-nickel concentrations.}
\label{empirical}
\end{figure}

\begin{figure}
\centering
\includegraphics[scale=0.36]{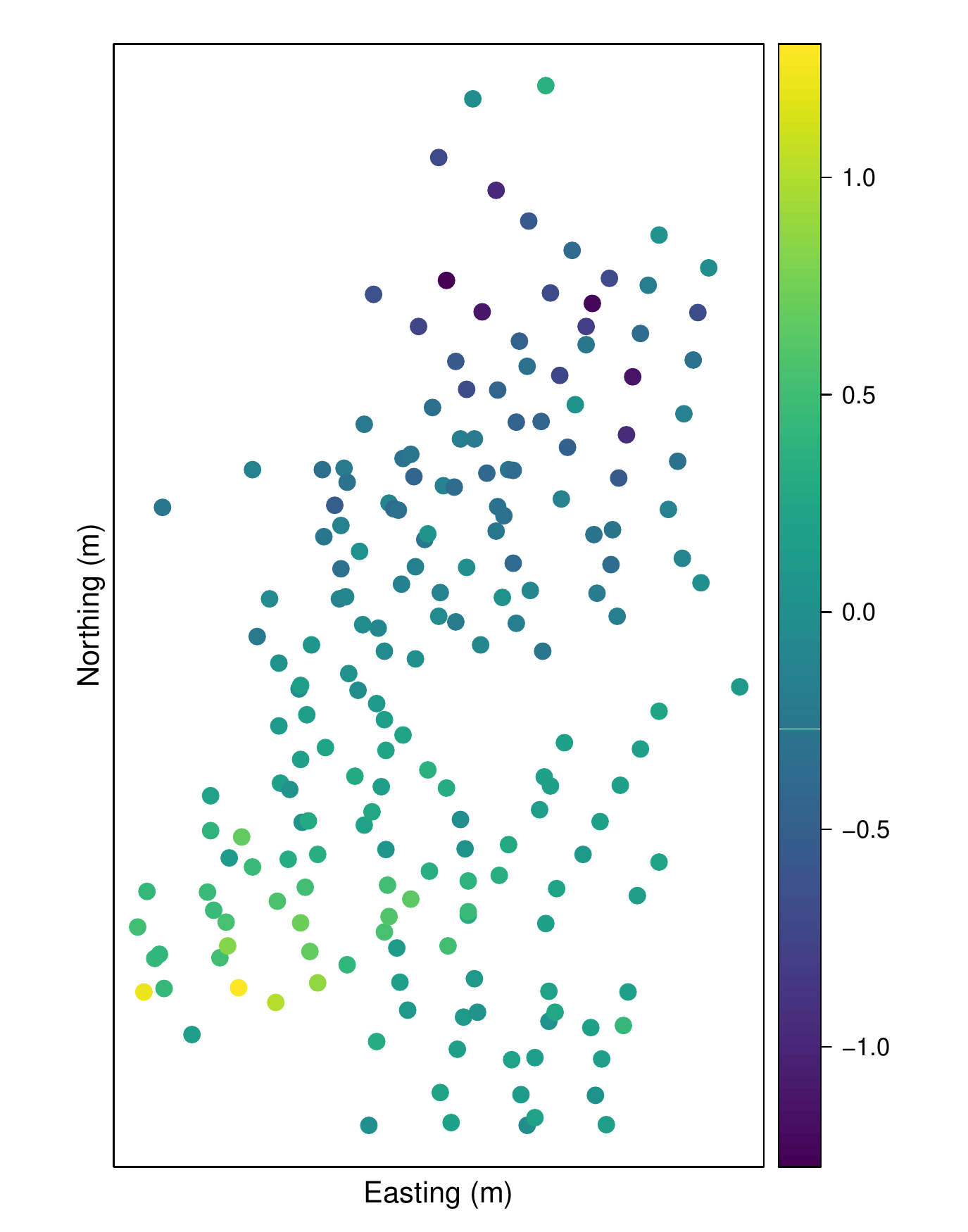}
\includegraphics[scale=0.375]{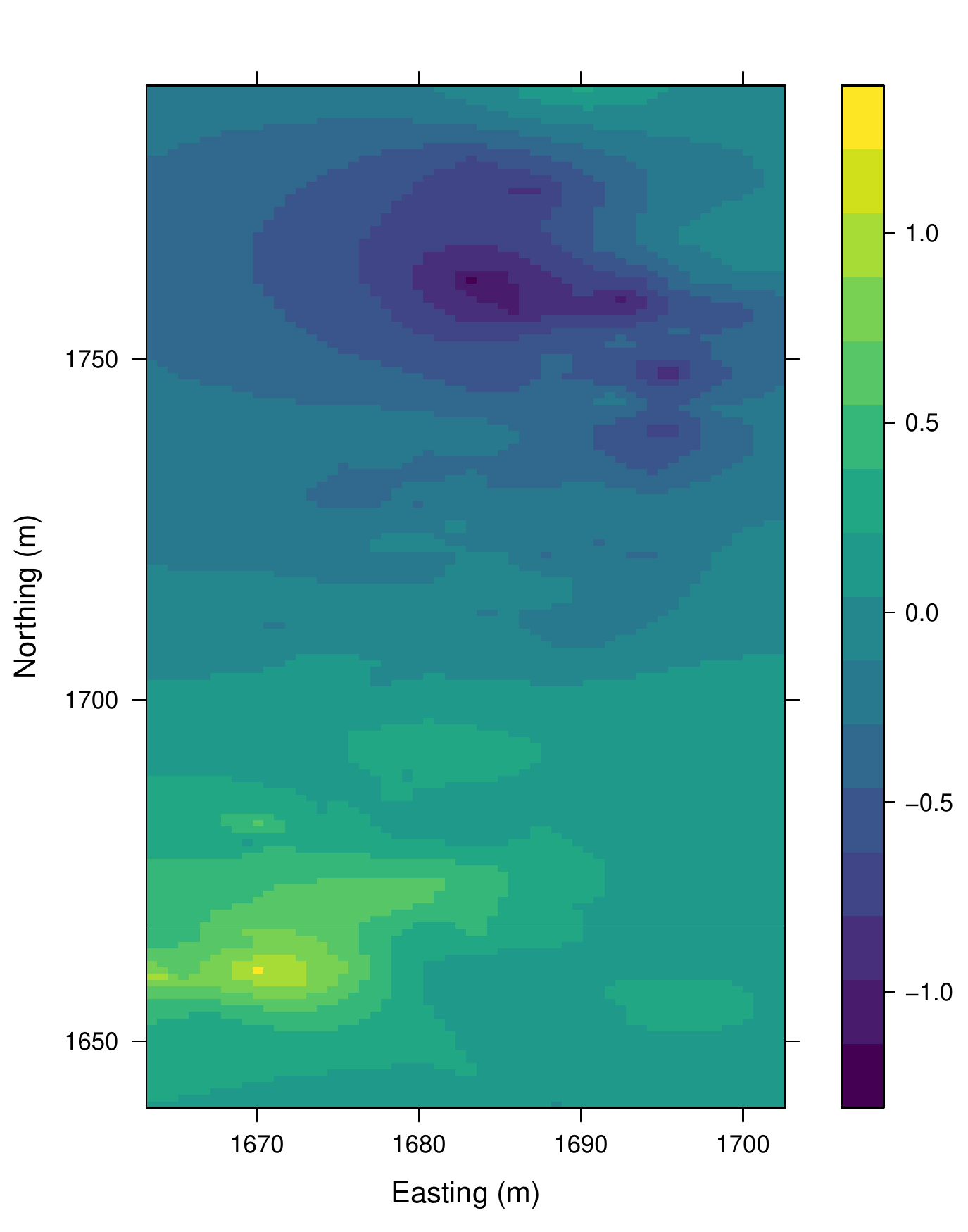}
\includegraphics[scale=0.375]{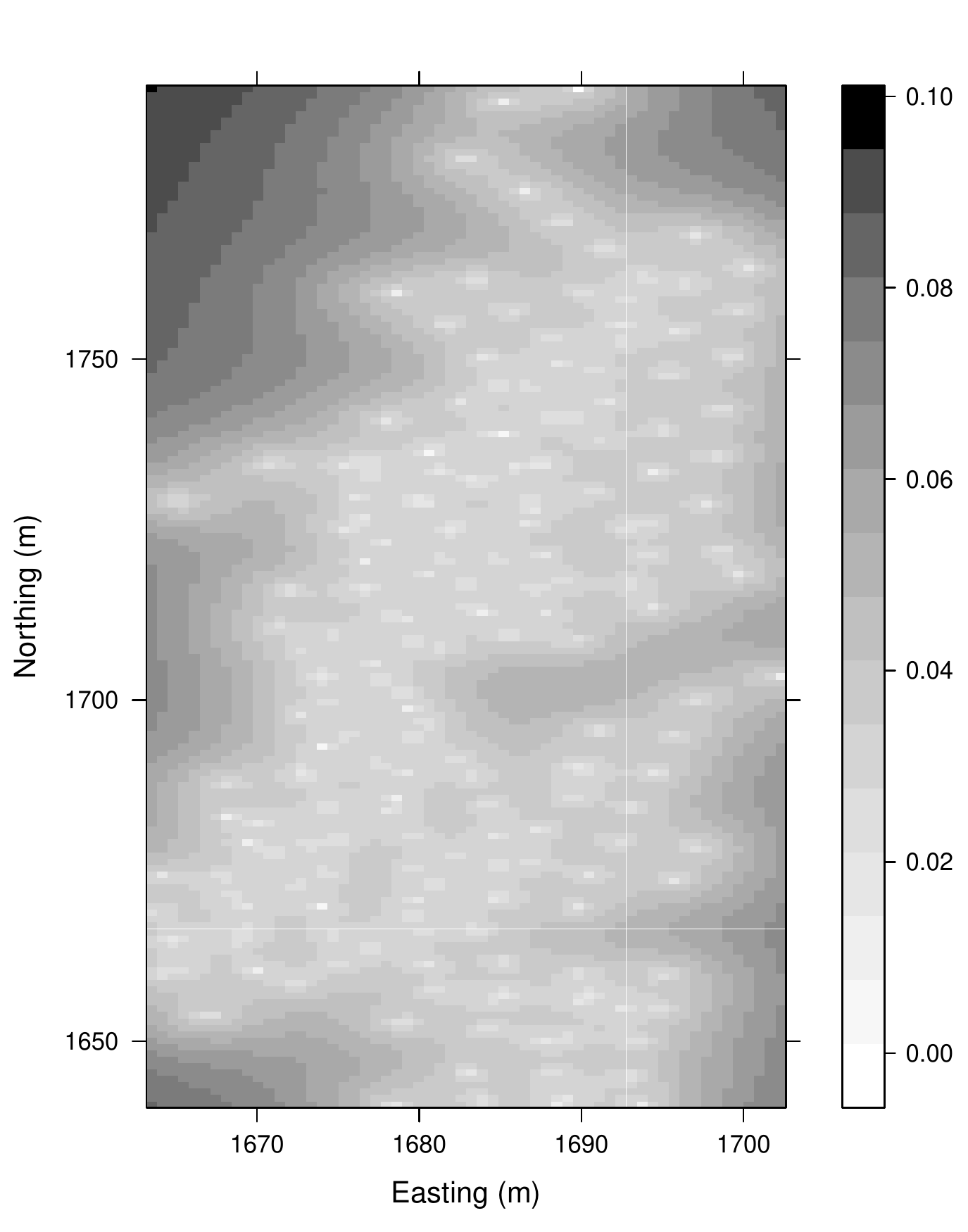}
\caption{Log-nickel concentrations (left), with the kriged surface (middle) and the corresponding variance (right).}
\label{heat_map}
\end{figure}

  This study focuses on nickel concentrations that are placed at an elevation of about $120$ meters, where $199$ irregularly spaced observations are available. We apply a log-transformation to reduce the skewness, and then the sample mean is subtracted. The resulting values are approximately Gaussian. The left panel of Figure \ref{heat_map} shows the transformed data set. 
   We fit two covariance models: the former is the parsimonious hybrid Cauchy-Matérn, parameterized by $\omega$, ${\alpha}$ and $\widetilde{\xi}$, with fixed $\nu_1=1/4$ and $\nu_2=1/2$, and the latter is the Generalized Cauchy, parameterized as in Section \ref{simula}, with fixed $\nu=1/4$ and $\delta=0.95$. The values of the fixed parameters have been selected after some experimental trials, taking into account the local behavior of the sample covariance (see Figure \ref{empirical}). 
 
 Table \ref{estimates_data} reports the likelihood estimates, with the corresponding standard errors, and the Akaike information criterion (AIC). We observe that the hybrid Cauchy-Matérn model outperforms its competitor in terms of AIC. Figure \ref{empirical} shows that the fitted covariance models  seem to be
reasonably close to the sample covariance. The fitted models differ substantially near the origin (distances less than 3 meters), since the hybrid model decays faster. On the contrary, for larger distances the hybrid model decays slower, although the difference between the curves becomes slight for distances greater than $15$ meters.

\begin{table}
\centering
\caption{Parameter estimates and Akaike Information Criterion (AIC) of fitted covariance models. Standard errors are reported in parentheses.}
\renewcommand{\arraystretch}{1.3}
\begin{tabular}{ccccc}  \hline \hline  
 Model & $\omega$ & $\alpha$ &  $\widetilde{\xi}$ &  AIC\\ \hline 
Hybrid Cauchy-Matérn &  $1.055$ & $12.31$ &  $0.063$  & $-34.03$   \\
 &  $(0.2516)$ & $(2.801)$ &  $(0.042)$  & \\
Generalized Cauchy &  $0.164$  & $1.959$  &  $-$   &  $-31.93$  \\
 &  $(0.040)$ & $(0.849)$ &  $-$  &  \\
\hline
\end{tabular}
\label{estimates_data}
\end{table}

In order to compare the models in terms of predictive performance, we conduct a cross-validation study, in a similar fashion to the experiments performed with simulated data.  Table \ref{scores_data} shows evidence, based on a leave-one-out cross-validation scheme, that the hybrid model has a better performance for this specific data set. In percentage terms, the MSE shows an improvement of approximately $3.4\%$. The largest difference occurs when we compare the LSCORE's (about $9\%$ improvement). 

We conclude this section with an illustration of a downscaled map of log-nickel concentrations (see Figure \ref{heat_map}), using the hybrid Cauchy-Matérn model. The interpolated spatial map, which is obtained through simple kriging,  is exhibited on a spatial grid of approximately $1$ meter ($7500$ locations). This kriged surface could be useful in small-scale mining processes, as it is a crucial step for industrial exploration and to quantify mineral reserves.

\begin{table}
\renewcommand{\arraystretch}{1.3}
\centering
\caption{Scores for the leave-one-out cross-validation study of log-nickel concentrations.}
\begin{tabular}{ccccc} \hline \hline
Model & MSE & MAE & LSCORE & CRPS\\ \hline 
Hybrid Cauchy-Matérn &    $0.0428$ &  $0.1431$ & $-0.1840$ & $0.4113$  \\
Generalized Cauchy &  $0.0443$  & $0.1462$  & $-0.1677$  & $0.4159$  \\
\hline
\end{tabular}
\label{scores_data}
\end{table}

\section{Conclusions and Perspectives} \label{sec5}
\label{conclusions}

We introduced a simple formalism to build sophisticated parametric families of  covariance functions.  We focused on a combination between the Matérn and Cauchy models, where local (mean square differentiability) and global (long memory) properties coexist in a single family.  We have also illustrated the use of our methodology by constructing a model that behaves as the Matérn class at short distances and attains negative values at large distances. Simulation studies show that a parsimonious hybrid Cauchy-Matérn model has statistically identifiable parameters. Also, this model provides improvements in terms of predictive performance in comparison to existing models, when a hybrid inherent dependence structure is present.  We reach similar conclusions when we apply this methodology to a mining dataset. While similar numerical studies could be performed for the hybrid Hole-Effect-Matérn model, we avoid them for the sake of simplicity and brevity. Additional interesting extensions of this work can be tackled in future investigations. We now provide two concrete research lines that could emerge from this work.

\subsection*{Multivariate Hybrid Covariance Models} In many practical situations, two or more variables are simultaneously recorded. Thus, our findings can be generalized to the case of multivariate fields $\{ {\bf Z}({\bf s})= (Z_1({\bf s}),\ldots, Z_p({\bf s}))^{\top}, \; {\bf s} \in \mathbb{R}^d \}$, having an isotropic matrix-valued covariance function $\Phi: [0,\infty) \to \mathbb{R}^{p \times p}$, that is,
$ {\rm cov}[Z_i({\bf s}),Z_j({\bf s}')] = \Phi_{ij}(h)$,  $h\geq 0$, where $h=\|{\bf s} - {\bf s}'\|$ and $i,j=1,\hdots,p$. We propose the hybrid model
\begin{equation*}
\label{expresion1}
   \widetilde{\Phi}(h; \bm{\lambda},\bm{\omega},\bm{\xi})   =  \omega_1   \int_{0}^{\xi_1}   \exp(-uh^2) {\bf G}_1(u;\bm{\lambda}_1)  \text{d}u + \omega_2   \int_{\xi_2}^{\infty}   \exp(-uh^2) {\bf G}_2(u;\bm{\lambda}_2)  \text{d}u,
\end{equation*}
that generalizes (\ref{expresion10}), where the vectors of parameters $\bm{\lambda}_i$ must be chosen in such a way that  the $p\times p$ matrices ${\bf G}_i(u;\bm{\lambda}_i)$ are positive semi-definite for every fixed $u \ge 0$. Hence, a straight application of Proposition 4 in \cite{porcu2011characterization} would ensure $\widetilde{\Phi}$ to be positive semi-definite. A multivariate version of the  hybrid Cauchy-Matérn covariance function is a natural candidate. The works of \cite{Gneiting:Kleibler:Schlather:2010} and \cite{moreva2022bivariate} are relevant to tackle this challenge. A multivariate version of the formulation (\ref{expresion10_gen}) could be deduced similarly. \medskip

\subsection*{Hybrid Covariance Models on Spheres} For random fields that are indexed by the $d$-dimensional unit sphere, $\mathbb{S}^d$, which is a useful framework when analyzing global data  ($\mathbb{S}^2$ is used as an approximation of the Earth), the isotropy assumption is given by ${\rm cov}[Z({\bf s}),Z({\bf s}')] = \psi(\theta)$, ${\bf s},{\bf s}'\in\mathbb{S}^d$, where $\psi:[0,\pi]\rightarrow \mathbb{R}$ is a continuous mapping and $\theta = \arccos({\bf s}^\top {\bf s}') \in [0,\pi]$ is the geodesic distance. Schoenberg's characterization \citep{schoenberg1942positive} establishes that a parametric isotropic covariance function $\psi(;\bm{\lambda})$ is valid in any dimension $d$, if and only if, it can be written as  ${\psi}(\theta; \bm{\lambda})  = \sum_{\ell=0}^\infty \beta_{\ell}(\bm{\lambda}) (\cos \theta)^{\ell}$,  $\theta\in[0,\pi]$, for some nonnegative and summable parametric sequence $\{\beta_{\ell}(\bm{\lambda}) \}_{\ell=0}^\infty$. Thus, the hybrid models can be adapted to the spherical context by considering a modified sequence of the form
 $$\widetilde{\beta}_{\ell}(\bm{\lambda},\bm{\omega},\bm{\xi})   =    \omega_1   \, \beta_{\ell}^{(1)}(\bm{\lambda}_1) 1_{ [0,  \lfloor \xi_1 \rfloor ) }(\ell)  +  \omega_2  \, \beta_{\ell}^{(2)}(\bm{\lambda}_2) 1_{ [ \lfloor \xi_{2}  \rfloor ,  \infty) }(\ell), \qquad \ell=0,1,\hdots,$$
where $ \lfloor \xi_{i} \rfloor \geq 0$, for $i=1,2$, with $\lfloor \cdot \rfloor$ standing for the floor function,  and $\beta_{\ell}^{(i)}$ being a nonnegative and summable sequence. The local properties of spherically indexed random fields, and their connections with the covariance function, have been studied in past literature \citep{bingham,guinness2016isotropic}. However, global properties such as long memory are less intuitive in this scenario as the spatial domain is a compact set. Covariance functions with hole effect, for low-dimensional spheres, could be obtained  by adapting formulation (\ref{expresion10_gen}).


\section*{Acknowledgements}

Alfredo Alegr\'ia was partially supported by the National Agency for Research and Development of Chile, through grant ANID/FONDECYT/INICIACIÓN/No.\,11190686. Fabián Ramirez was partially supported by the Direcci\'on de Postgrados y Programas (DPP) of the Universidad Técnica Federico
Santa María.  Emilio Porcu is supported by the Khalifa University of Science and Technology under Award No. FSU-2021-016.


\bibliographystyle{apalike}
\bibliography{mybib}

\begin{thebibliography}{}

\bibitem[Abramowitz and Stegun, 1972]{Abramowitz-Stegun:1965}
Abramowitz, M. and Stegun, I.~A. (1972).
\newblock {\em Handbook of Mathematical Functions with Formulas, Graphs, and
  Mathematical Tables}.
\newblock Dover Publications.

\bibitem[Alegr{\'\i}a, 2020]{alegria2020crossdimple}
Alegr{\'\i}a, A. (2020).
\newblock Cross-dimple in the cross-covariance functions of bivariate isotropic
  random fields on spheres.
\newblock {\em Stat}, 9(1):e301.

\bibitem[Alegr{\'\i}a et~al., 2021]{alegria2021bivariate}
Alegr{\'\i}a, A., Emery, X., and Porcu, E. (2021).
\newblock Bivariate {M}at{\'e}rn covariances with cross-dimple for modeling
  coregionalized variables.
\newblock {\em Spatial Statistics}, 41:100491.

\bibitem[Barndorff-Nielsen, 1978]{barndorff1978hyperbolic}
Barndorff-Nielsen, O. (1978).
\newblock Hyperbolic distributions and distributions on hyperbolae.
\newblock {\em Scandinavian Journal of Statistics}, 5(3):151--157.

\bibitem[Barp et~al., 2022]{barp2022riemann}
Barp, A., Oates, C.~J., Porcu, E., and Girolami, M. (2022).
\newblock A {R}iemann--{S}tein kernel method.
\newblock {\em Bernoulli}, 28(4):2181--2208.

\bibitem[Berg et~al., 2008]{berg2008dagum}
Berg, C., Mateu, J., and Porcu, E. (2008).
\newblock The {D}agum family of isotropic correlation functions.
\newblock {\em Bernoulli}, 14(4):1134--1149.

\bibitem[Bevilacqua et~al., 2022]{bevilacqua2022unifying}
Bevilacqua, M., Caama{\~n}o-Carrillo, C., and Porcu, E. (2022).
\newblock Unifying compactly supported and {M}at{\'e}rn covariance functions in
  spatial statistics.
\newblock {\em Journal of Multivariate Analysis}, 189:104949.

\bibitem[Bevilacqua and Gaetan, 2015]{bevilacqua2015comparing}
Bevilacqua, M. and Gaetan, C. (2015).
\newblock Comparing composite likelihood methods based on pairs for spatial
  {G}aussian random fields.
\newblock {\em Statistics and Computing}, 25(5):877--892.

\bibitem[Bingham, 1973]{bingham}
Bingham, N.~H. (1973).
\newblock Positive {D}efinite {F}unctions on {S}pheres.
\newblock {\em Proc. Cambridge Phil. Soc.}, 73:145--156.

\bibitem[Chaudhry and Zubair, 1994]{chaudhry1994generalized}
Chaudhry, M.~A. and Zubair, S.~M. (1994).
\newblock Generalized incomplete gamma functions with applications.
\newblock {\em Journal of Computational and Applied Mathematics},
  55(1):99--123.

\bibitem[Chen et~al., 2018]{chen2018current}
Chen, W., Castruccio, S., Genton, M.~G., and Crippa, P. (2018).
\newblock Current and future estimates of wind energy potential over {S}audi
  {A}rabia.
\newblock {\em Journal of Geophysical Research: Atmospheres},
  123(12):6443--6459.

\bibitem[Chilès and Delfiner, 2012]{Chiles2012}
Chilès, J.-P. and Delfiner, P. (2012).
\newblock {\em Geostatistics: {M}odeling {S}patial {U}ncertainty, {S}econd
  {E}dition}.
\newblock John Wiley \& Sons.

\bibitem[Cockayne et~al., 2019]{cockayne2019bayesian}
Cockayne, J., Oates, C.~J., Sullivan, T.~J., and Girolami, M. (2019).
\newblock Bayesian probabilistic numerical methods.
\newblock {\em SIAM review}, 61(4):756--789.

\bibitem[Cressie, 1993]{Cressie:1993}
Cressie, N. (1993).
\newblock {\em Statistics for Spatial Data}.
\newblock Wiley, New York, revised edition.

\bibitem[Cressie and Kornak, 2003]{cressie2003spatial}
Cressie, N. and Kornak, J. (2003).
\newblock Spatial {S}tatistics in the {P}resence of {L}ocation {E}rror with an
  {A}pplication to {R}emote {S}ensing of the {E}nvironment.
\newblock {\em Statistical Science}, 18(4):436--456.

\bibitem[Daley and Porcu, 2014]{daley2014dimension}
Daley, D. and Porcu, E. (2014).
\newblock Dimension walks and {S}choenberg spectral measures.
\newblock {\em Proceedings of the American Mathematical Society},
  142(5):1813--1824.

\bibitem[Di~Lorenzo et~al., 2014]{horma2}
Di~Lorenzo, E., Combes, V., Keister, J., Strub, P., Andrew, T., Peter, F.,
  Marck, O., Furtado, J.~., Bracco, A., Bograd, S., Peterson, W., Schwing, F.,
  Taguchi, B., Hormázabal, S., and Parada, C. (2014).
\newblock Synthesis of {P}acific {O}cean {C}limate and ecosystem dynamics.
\newblock {\em Oceanography}, 26(4):68--81.

\bibitem[Edwards et~al., 2019]{sacowea1}
Edwards, M., Castruccio, S., and Hammerling, D. (2019).
\newblock A multivariate {G}lobal {S}patio-{T}emporal {S}tochastic {G}enerator
  for {C}limate {E}nsembles.
\newblock {\em Journal of Agricultural, Biological and Environmental Sciences}.

\bibitem[Emery and Lantu{\'e}joul, 2006]{emery2006tbsim}
Emery, X. and Lantu{\'e}joul, C. (2006).
\newblock Tbsim: A computer program for conditional simulation of
  three-dimensional gaussian random fields via the turning bands method.
\newblock {\em Computers \& Geosciences}, 32(10):1615--1628.

\bibitem[Emery and S{\'e}guret, 2020]{emery2020geostatistics}
Emery, X. and S{\'e}guret, S.~A. (2020).
\newblock {\em Geostatistics for the Mining Industry: Applications to Porphyry
  Copper Deposits}.
\newblock CRC Press.

\bibitem[Furrer et~al., 2006]{furrer2006covariance}
Furrer, R., Genton, M.~G., and Nychka, D. (2006).
\newblock Covariance tapering for interpolation of large spatial datasets.
\newblock {\em Journal of Computational and Graphical Statistics},
  15(3):502--523.

\bibitem[Furrer et~al., 2007]{furrer2007multivariate}
Furrer, R., Sain, S.~R., Nychka, D., and Meehl, G.~A. (2007).
\newblock Multivariate bayesian analysis of atmosphere--ocean general
  circulation models.
\newblock {\em Environmental and Ecological Statistics}, 14(3):249--266.

\bibitem[Gneiting et~al., 2010]{Gneiting:Kleibler:Schlather:2010}
Gneiting, T., Kleiber, W., and Schlather, M. (2010).
\newblock Mat\'ern {cross-covariance} functions for multivariate random fields.
\newblock {\em Journal of the American Statistical Association},
  105:1167--1177.

\bibitem[Gneiting and Schlather, 2004]{gneiting2004stochastic}
Gneiting, T. and Schlather, M. (2004).
\newblock Stochastic models that separate fractal dimension and the {H}urst
  effect.
\newblock {\em SIAM review}, 46(2):269--282.

\bibitem[Gradshteyn and Ryzhik, 2007]{Grad}
Gradshteyn, I. and Ryzhik, I. (2007).
\newblock {\em Table of {I}ntegrals, {S}eries, and {P}roducts}.
\newblock Academic Press, Amsterdam, 7th. edition.

\bibitem[Guinness and Fuentes, 2016]{guinness2016isotropic}
Guinness, J. and Fuentes, M. (2016).
\newblock Isotropic covariance functions on spheres: {S}ome properties and
  modeling considerations.
\newblock {\em Journal of Multivariate Analysis}, 143:143--152.

\bibitem[Guinness and Hammerling, 2018]{guinness2018compression}
Guinness, J. and Hammerling, D. (2018).
\newblock Compression and conditional emulation of climate model output.
\newblock {\em Journal of the American Statistical Association},
  113(521):56--67.

\bibitem[Hristopulos, 2020]{hristopulos2020random}
Hristopulos, D.~T. (2020).
\newblock {\em Random fields for spatial data modeling}.
\newblock Springer.

\bibitem[James et~al., 2013]{james2013introduction}
James, G., Witten, D., Hastie, T., and Tibshirani, R. (2013).
\newblock {\em An introduction to statistical learning}, volume 112.
\newblock Springer.

\bibitem[Kaufman et~al., 2008]{kaufman2008covariance}
Kaufman, C.~G., Schervish, M.~J., and Nychka, D.~W. (2008).
\newblock Covariance tapering for likelihood-based estimation in large spatial
  data sets.
\newblock {\em Journal of the American Statistical Association},
  103(484):1545--1555.

\bibitem[Laga and Kleiber, 2017]{laga2017modified}
Laga, I. and Kleiber, W. (2017).
\newblock The modified {M}at{\'e}rn process.
\newblock {\em Stat}, 6(1):241--247.

\bibitem[Ma and Bhadra, 2022]{ma2022beyond}
Ma, P. and Bhadra, A. (2022).
\newblock Beyond {M}at{\'e}rn: {O}n a {C}lass of {I}nterpretable {C}onfluent
  {H}ypergeometric {C}ovariance {F}unctions.
\newblock {\em Journal of the American Statistical Association}, (in press).

\bibitem[Mat\'ern, 1986]{Matern}
Mat\'ern, B. (1986).
\newblock {\em Spatial {V}ariation — {S}tochastic {M}odels and {T}heir
  {A}pplication to {S}ome {P}roblems in {F}orest {S}urveys and {O}ther
  {S}ampling {I}nvestigations}.
\newblock Springer.

\bibitem[Moreva and Schlather, 2022]{moreva2022bivariate}
Moreva, O. and Schlather, M. (2022).
\newblock Bivariate covariance functions of {P}{\'o}lya type.
\newblock {\em Journal of Multivariate Analysis}, in press.

\bibitem[Ostoja-Starzewski, 2006]{ostoja2006material}
Ostoja-Starzewski, M. (2006).
\newblock Material spatial randomness: From statistical to representative
  volume element.
\newblock {\em Probabilistic engineering mechanics}, 21(2):112--132.

\bibitem[Pazouki and Schaback, 2011]{pazouki2011bases}
Pazouki, M. and Schaback, R. (2011).
\newblock Bases for kernel-based spaces.
\newblock {\em Journal of Computational and Applied Mathematics},
  236(4):575--588.

\bibitem[Porcu et~al., 2018a]{porcu2018modeling}
Porcu, E., Alegria, A., and Furrer, R. (2018a).
\newblock Modeling temporally evolving and spatially globally dependent data.
\newblock {\em International Statistical Review}, 86(2):344--377.

\bibitem[Porcu et~al., 2018b]{porcu2018shkarofsky}
Porcu, E., Bevilacqua, M., and Hering, A.~S. (2018b).
\newblock The {S}hkarofsky-{G}neiting class of covariance models for bivariate
  {G}aussian random fields.
\newblock {\em Stat}, 7(1):e207.

\bibitem[Porcu and Zastavnyi, 2011]{porcu2011characterization}
Porcu, E. and Zastavnyi, V. (2011).
\newblock Characterization theorems for some classes of covariance functions
  associated to vector valued random fields.
\newblock {\em Journal of Multivariate Analysis}, 102(9):1293--1301.

\bibitem[Posa, 2022]{posa}
Posa, D. (2022).
\newblock Special classes of isotropic covariance functions.
\newblock {\em Preprint, Research Square}.

\bibitem[Schaback and Wendland, 2006]{schaback2006kernel}
Schaback, R. and Wendland, H. (2006).
\newblock Kernel techniques: from machine learning to meshless methods.
\newblock {\em Acta numerica}, 15:543--639.

\bibitem[Schlather, 2010]{schlather2010some}
Schlather, M. (2010).
\newblock Some covariance models based on normal scale mixtures.
\newblock {\em Bernoulli}, 16(3):780--797.

\bibitem[Schlather and Moreva, 2017]{schlather2017parametric}
Schlather, M. and Moreva, O. (2017).
\newblock A parametric model bridging between bounded and unbounded variograms.
\newblock {\em Stat}, 6(1):47--52.

\bibitem[Schoenberg, 1938]{schoenberg1938}
Schoenberg, I.~J. (1938).
\newblock Metric {S}paces and {C}ompletely {M}onotone {F}unctions.
\newblock {\em Annals of Mathematics}, 39(4):811--841.

\bibitem[Schoenberg, 1942]{schoenberg1942positive}
Schoenberg, I.~J. (1942).
\newblock Positive definite functions on spheres.
\newblock {\em Duke Mathematical Journal}, 9(1):96--108.

\bibitem[Stein, 1999]{stein-book}
Stein, M.~L. (1999).
\newblock {\em Statistical {I}nterpolation of {S}patial {D}ata: {S}ome {T}heory
  for {K}riging}.
\newblock Springer, New York.

\bibitem[Stein, 2007]{stein2007spatial}
Stein, M.~L. (2007).
\newblock Spatial variation of total column ozone on a global scale.
\newblock {\em The Annals of Applied Statistics}, 1(1):191--210.

\bibitem[Yaglom, 1987]{yaglom1987correlation}
Yaglom, A.~M. (1987).
\newblock {\em Correlation Theory of Stationary and Related Random Functions,
  Volume I: Basic Results}, volume 131.
\newblock Springer.

\bibitem[Zhang and Wang, 2010]{zhang2010kriging}
Zhang, H. and Wang, Y. (2010).
\newblock Kriging and cross-validation for massive spatial data.
\newblock {\em Environmetrics}, 21(3-4):290--304.

\end{thebibliography}

\end{document}